\def\qed{\hfill{$\diamond$}}
\def\v(#1,#2,#3,#4){\pmatrix{#1\cr #2\cr #3\cr #4}}
\def\pf(#1,#2){\langle\langle #1,#2\rangle\rangle}
\def\pf(#1){\langle\langle #1\rangle\rangle}
\def\symb(#1){\langle #1\rangle}
\def\Hom{\mathop{\rm Hom}}
\def\sgn{\mathop{\rm sgn}}
\def\inv{\mathop{\rm inv}}
\def\supp{\mathop{\rm supp}}

\def\Z{{\bf Z}}
\def\R{{\bf R}}
\def\P{{\bf P}}

\def\cC{{C_*}}
\def\cK{{K_*}}
\def\cP{{P_*}}
\def\dotK{\dot{K}}
\def\sqcl{{\dotK/\dotK^2}}

\def\Zm(#1){{\bf Z}/#1}
\def\Ext(#1,#2){\mathop{\rm Ext}(#1,#2)}
\def\Q{{\bf Q}}

\def\f{{\bf f}}

\def\broose{\atopwithdelims[]}
\def\qed{\hfill{$\diamond$}}
\def\v(#1,#2,#3,#4){\pmatrix{#1\cr #2\cr #3\cr #4}}
\def\pf(#1,#2){\langle\langle #1,#2\rangle\rangle}
\def\pf(#1){\langle\langle #1\rangle\rangle}
\def\symb(#1){\langle #1\rangle}
\def\Hom{\mathop{\rm Hom}}
\def\sgn{\mathop{\rm sgn}}
\def\inv{\mathop{\rm inv}}
\def\supp{\mathop{\rm supp}}

\def\Z{{\bf Z}}
\def\R{{\bf R}}
\def\P{{\bf P}}
\def\Pp{{\bf P}_+}

\def\cC{{C_*}}
\def\cK{{K_*}}
\def\cP{{P_*}}
\def\dotK{\dot{K}}
\def\sqcl{{\dotK/\dotK^2}}

\def\Zm(#1){{\bf Z}/#1}
\def\Ext(#1,#2){\mathop{\rm Ext}(#1,#2)}
\def\Q{{\bf Q}}

\def\f{{\bf f}}

\def\G{{G}}
\def\Gp{{G}_+}
\def\PG{{PG}}
\def\PpG{{P_+G}}
\def\PGp{{PG}_+}
\def\PpGp{{P_+G_+}}
\def\eu{{eu}}
\newcount\m
\def\t{\the\m\global\advance\m by 1}
\m=1
\newcount\n
\def\f{\the\n\global\advance\n by 1}
\n=1

\centerline{\bf Tautological characteristic classes I.}
\bigskip
\rm
Jan Dymara

Instytut Matematyczny Uniwersytetu Wroc\l awskiego,

pl.~Grunwaldzki 2, 50-384 Wroc\l aw
\smallskip

Tadeusz Januszkiewicz

Instytut Matematyczny PAN,

ul.~Kopernika 18, 51-617 Wroc\l aw

\footnote{}{Both authors were supported by
Polish NCN grant UMO-2016/23/B/ST1/01556.}

\bigskip\rm

\centerline{\bf Introduction.}
\medskip
Take a chain complex $\cC$  and fix the degree $k$.
The identity map $C_k \to C_k$ can be viewed as a cochain
of degree $k$, with values/coefficients in $C_k$.
Usually it is not a cocycle, but we can force it to be,
dividing by boundaries: changing coefficients to $C_k/B_k$.
Denote the resulting cocycle by $T$.

Now suppose that a group $G$ acts on $\cC$. Clearly $T$ is $G$-equivariant,
in other words it is a cocycle with twisted coefficients.
We can force $T$  to be a constant coefficients cocycle simply
dividing the coefficients further down to 
the biggest quotient of the $G$-module $C_k/B_k$ on which $G$ acts trivially,
called coinvariants of $G$.
Denote the resulting image of $T$ by $\tau$.

Besides producing an untwisted cocycle, this construction has an additional 
crucial advantage:
the modules $C_k$, or even $C_k/B_k$, are usually very big, while the coinvariants
$(C_k/B_k)_G$ are much smaller and sometimes manageable.

It is of interest to go half-way in this procedure:
fix a (large) normal subgroup $N$ of $G$, and take coinvariants $(C_k/B_k)_N$.
Then $T$ becomes a (slightly twisted by an action of $G/N$) cocycle 
taking values in (sometime still manageable, but bigger) module of $N$-coinvariants.

One has every reason to expect that this purely algebraic construction 
has nice functorial properties, and that it carries a significant amount
of information about $\cC$ as a $G$-module.
Theorem 1.5 spells out the most natural form of functoriality.

This algebraic construction needs an input.
For us such an input comes from a geometry (or, as some will undoubtedly insist, algebra),
namely from the complex of geometric configurations.
One takes a homogeneous space $G/H$ (for example a projective space over an arbitrary field $K$)
and builds a simplicial complex whose simplices are  $n$-tuples of points ``in general position". 
The notion of general position that we use depends on the situation, and is discussed separately in each case,
but the underlying idea is uniform.
In all cases the simplicial complexes that we consider have an additional crucial {\it ``star property"}
which is discussed in Section 2. The star property makes the simplicial complex contractible
in a strong, geometric sense. 
The chain complex is just the complex of (alternating or ordered) simplicial chains.

In order to construct characteristic classes of flat $G$-bundles, 
we have to address the problem that the $G$-action on the space of configurations is not free.
But this is done in a standard way, by ``Borel construction", that we execute
on the chain level.
We end up with cocycles living in a cochain complex computing group cohomology
of $G$ (seen as a discrete group). 
The star property implies almost immediately that the cohomology class of the cocycle is bounded (cf.~Theorem 3.1).
The boundedness of tautological
classes is taken with respect to the natural seminorm on the coefficient group.
One should keep in mind that this has fairly different overtones from the usual $\R$-coefficients
bounded cohomology of [Gro], compare for example [Ghys]. 

The first interesting case when this construction produces something valuable is 
that of $PGL(2,K)$ acting on the projective line. 
This has been studied by Nekov\'a$\check{\rm r}$, who defined and studied the ``Witt class'' for $PSL(2,K)$, 
with coefficients in the Witt ring of quadratic forms over $K$. 
It is an amazing fact that the four term Witt relation $[a]+[b]=[a+b] +[ab(a+b)]$ is indeed the cocycle
relation in the complex of projective point configurations. We review this in detail in Section 7.

The main results of the present paper concern the construction and study of the
``Euler class for flat $PGL(n,K)$-bundles" in the case where
$K$ is an arbitrary ordered field and $n$ is even.
This class is constructed using the general strategy outlined above.
We take the $PGL(n,K)$ action on the simplicial
complex of generic configurations of points in $\P^{n-1}(K)$,
the induced action on $C_n/B_n$, and then we take coinvariants
with respect to the group $PGL_+(n,K)$ of maps with {\it positive} determinant.
(Note that coinvariants with respect to the full projective group are trivial,
while coinvariants with respect to $PSL(n,K)$ are too large for us to
handle---for $PGL_+(n,K)$ we have a nice answer.)
The resulting tautological class $\eu$
is (an analogue of) the Euler class---for  flat $PGL(n,K)$-bundles. 
It is twisted by the homomorphism to $\Z/2$ whose kernel consists of maps
with positive determinant. The coefficients are $\Z$ for $n$ even and
trivial for
$n$ odd (cf.~Theorem 8.1).

One can run a parallel construction starting from
the $GL_+(n,K)$ action on the {\it positive} projective space $\P_+^{n-1}(K)$.
The resulting class $\eu_+$ has coefficients in a free abelian group
of rank $\lfloor n/2\rfloor+1$ (cf.~Theorem 8.1; admittedly, the computation
here is somewhat heavy). Consequently, $\eu_+$ can be split
into components $\eu_k$ that are cohomology classes with $\Z$ coefficients.

We prove several results about the Euler classes $\eu$ and $\eu_+$.
Theorems 9.1 and 10.1 explain the relation between various
components of $eu_+$.
Theorem 11.1 gives a clean formula for the Euler class of a
cross product of bundles, while Theorem 12.5
gives a cup product formula for the direct sum.
In Section 13 we discuss functoriality. In particular, we relate
$\eu$ and $\eu_+$ in Theorem 13.1.
We also compare the
Euler and Witt classes for $PSL(2,K)$-bundles in Theorem 13.4. 
Finally, in Theorem 13.6, using the cross-product formula, we show
non-triviality of our Euler classes in every even dimension. 

Further characteristic classes and more applications are postponed to subsequent papers.

Tautological classes with coefficients
in $C_k/B_k$ were defined in a forgotten paper of James Dugundji [Dug],
where he also proved some results of general nature,
like functoriality and universality.
The paper was forgotten, probably because the results did not help with
actual calculations:
modules $C_k/B_k$ are usually very big and unmanageable.
We discovered Dugundji's paper when we were already well into our project.
Our initial inspiration came from the papers of
Nekov\'a$\check{\rm r}$ [Ne] and Kramer and Tent [Kr-T], where
the idea of passing to $G$-coinvariants is present. With a grain of salt,
one may say that the Witt and Maslov classes are constructed
in these papers in the tautological way.

Reznikov ([Rez]) noticed that for an ordered field $K$
one has an ``Euler class'' for $PSL(2,K)$ with $\Z$ coefficients.
In fact, this class is (a multiple of)
the image of the Witt class of Nekov\'a$\check{\rm r}$
under the signature map from the Witt ring to $\Z$,
given by the ordering of $K$.

\medskip
The plan of the paper is as follows.

In Part I we discuss the general theory:
Section 1 explains definitions and functoriality of tautological classes in a
purely algebraic, abstract context; in Section 2 it is shown how actions on simplicial
complexes can lead to examples, star-property is recalled, and a method of coefficient calculation
for actions on simplicial complexes is described; Section 3 is about
(automatic) boundedness of tautological classes; Section 4 contains a simplicial counterpart of
the process of representing classes of flat bundles by pull-backs of invariant forms via sections.

Part II is about $GL(2)$: in Section 6 we discuss various actions of this group
with a view towards investigating the corresponding tautological classes;
in Section 7 the Witt group appears as the coefficient group coming from
the general formalism applied to the homographic action on the projective line,
and the tautological Witt class is defined.

In Part III we define Euler classes for the groups $PGL(n,K)$ and $PGL_+(n,K)$,
where $K$ is an arbitrary ordered field. In Section 8 actions of these groups
on $\P^{n-1}(K)$ and on $\P^{n-1}_+(K)$ are used to define tautological
Euler classes $\eu$ and $\eu_+$; coefficients are calculated, and $\eu_+$
is decomposed into a direct sum of classes $\eu_k$ (with coefficient in $\Z$).
In Section 9 we establish a general relation between the classes $\eu_k$,
and in Section 10 we express all of them in terms of $\eu_0$ in a weak sense
using Smillie's argument. In Sections 11 and 12 we show some
multiplicativity properties of $\eu_0$. In Section 13 we further investigate
relations between various Euler classes (and the Witt class); we also prove
that all these classes are non-trivial (for $n$ even).
\medskip

We would like to acknowledge extensive discussions with Linus Kramer
that were very helpful in the initial stages of this project.

\rm

\bigskip
\vfill\eject
\centerline{\bf I. Generalities}
\bigskip
\bf 1. Algebraic tautological classes.\rm
\medskip
\def\a{1}
\m=1
\n=1
\bf Chain complexes. \rm Let $\cC=(C_n,\partial_n)$ be a chain complex of abelian groups.
As usual, we put $Z_n=\ker{\partial_n}$ (cycles) and
$B_n=\mathop{\rm im}{\partial_n}$ (boundaries).
Let us fix an integer $d$ and consider ${\rm id}_{C_d}$ as an element of
$\Hom(C_d,C_d)$---the $d$-cochain group of the complex $\Hom(\cC,C_d)$.
This element is usually not a cocycle---yet, if we replace the coefficient group
$C_d$ by the quotient $C_d/B_d$, it becomes one.
\smallskip
\bf Definition \a.\t.\rm\par
Let $\cC$ be a chain complex. The tautological cocycle $T_\cC^d$ is
the $d$-cycle of the complex $\Hom(\cC,C_d/B_d)$ defined by the quotient map
$C_d\to C_d/B_d$. The tautological class $\tau^d_\cC$ is the cohomology class of
$T^d_\cC$ in $H^d(\Hom(\cC,C_d/B_d))$.
\smallskip
\rm
The cochain $T^d_\cC$ is indeed a cocycle:
$$\delta T^d_\cC(c)=T^d_\cC(\partial c)=\partial c+B_d=B_d.$$

Notice that $\tau^d_\cC$ is functorial---in the following way:
Let $f\colon \cC\to\cK$ be a chain map. Then $f_d\colon C_d\to K_d$ induces
a map $C_d/\ker{\partial_d}\to K_d/\ker{\partial_d}$, which in turn induces
a map
$$f_*\colon H^d(\Hom(\cC,C_d/\ker{\partial_d}))\to H^d(\Hom(\cC,K_d/\ker{\partial_d})).$$
There is also the map $f^*\colon \Hom(\cK,K_d/\ker{\partial_d})\to \Hom(\cC,K_d/\ker{\partial_d})$
inducing
$$f^*\colon H^d(\Hom(\cK,K_d/\ker{\partial_d}))\to H^d(\Hom(\cC,K_d/\ker{\partial_d})).$$
Clearly, $f^*\tau^d_\cK=f_*\tau^d_\cC$:
indeed, 
both these classes are represented by the same cocycle
$C_c\ni c\mapsto f_d(c)+\ker{\partial_d}\in K_d/\ker{\partial_d}$.
\medskip
\bf $G$-chain complexes. \rm Now suppose that the complex $\cC$ is a $G$-chain complex, i.e., it is
acted upon by a group $G$, by chain maps. The group $C_d/B_d$ has the induced $G$-module structure.
The tautological
cocycle  $T^d_\cC\colon C_d\to C_d/B_d$ is a $G$-map.
\smallskip
\bf Definition \a.\t.\rm\par
Let $\cC$ be a $G$-chain complex. The tautological class
$\tau^d_{\cC,G}\in H^d(\Hom_G(\cC,C_d/B_d))$ (cohomology with twisted coefficients) is the cohomology class of $T^d_\cC$.
\smallskip
\rm
We have found out that the above class has also been defined and investigated in a forgotten paper of Dugundji [Dug].

The $G$-module $C_d/B_d$ is usually very big. To cut it down in size we
will consider its
coinvariants group $U_d=(C_d/B_d)_G$---its largest $G$-trivial quotient.
This group might be either too small to carry information or too big to extract information, yet in some cases
it is non-trivial and manageable.
\smallskip
\bf Definition \a.\t.\rm\par
Let $\cC$ be a $G$-chain complex. Let $U_d$ (or $U_d(\cC)$) denote the
coinvariants group $(C_d/B_d)_G$. The tautological class
$\tau^d_{\cC/G}\in H^d(\Hom_G(\cC,U_d))$ is the cohomology class
of $T^d_{\cC/G}$: the cocycle 
obtained by composing the tautological cocycle $T^d_\cC$ with the quotient map
$C_d/B_d\to (C_d/B_d)_G$.
\smallskip
\rm
\bf Remarks.\rm\par
\item{1)} {The functor of coinvariants is right-exact ([Brown, II, \S 2]).
Therefore we have $U_d=(C_d)_G/(B_d)_G$ (strictly speaking, we divide by the
image--not necessarily injective---of $(B_d)_G$ in $(C_d)_G$).
Moreover, $\partial\colon C_{d+1}\to C_d$ induces a map
$\partial\colon (C_{d+1})_G\to (C_d)_G$, and $U_d$ can also
be described as $(C_d)_G/\partial((C_{d+1})_G)$.}
\item{2)}
{If $N$ is a normal subgroup of $G$ then there exists yet another, $G/N$-twisted tautological class
$\tau$ in $H^d(\Hom_G(\cC,(C_d/B_d)_N))$.} 
\smallskip
Let us discuss functoriality. Suppose that $\cC$ is a $G$-complex and that $\cK$ is an $H$-complex.
Assume that $\phi\colon G\to H$ is a homomorphism and that $f\colon \cC\to \cK$ is a
$\phi$-equivariant chain map. The group $U_d(\cK)$ acquires a $G$-module structure via $\phi$.
We have two maps:
$$
H^d(\Hom{\!}_H(\cK,U_d(\cK)))
\mathop{\longrightarrow}\limits^{f^*}
H^d(\Hom{\!}_G(\cC,U_d(\cK)))
\mathop{\longleftarrow}\limits^{f_*}
H^d(\Hom{\!}_G(\cC,U_d(\cC))),
$$
the right one induced by the $f$-induced coefficient map $U_d(\cC)\to U_d(\cK)$.
As before, it is straightforward to check that
$f^*\tau^d_{\cK/H}=f_*\tau^d_{\cC/G}$---both these classes are represented
by the cocycle $C_d\ni c\mapsto [f_d(c)]\in U_d(\cK)$.
\medskip
\bf Acyclic $G$-chain complexes. \rm 
Let us now assume that $\cC$ is an acyclic $G$-chain complex. By this we mean that:
(1) $C_n=0$ for $n<0$; (2) $\cC$ comes equipped with an augmentation map---a $G$-homomorphism
$\epsilon\colon C_0\to\Z$, where $\Z$ has the trivial $G$-module structure; (3)
the augmented complex
$$\ldots\to C_n\to C_{n-1}\to\ldots\to C_1\to C_0\mathop{\to}\limits^{\epsilon}\Z\to 0\to\ldots$$
is exact. (In other words: $\cC$ is a resolution of the trivial $G$-module $\Z$.)

The tautological class $\tau^d_{\cC/G}$ can be used to define
a cohomology class of the group $G$, as follows.
Let $\cP$ be a projective resolution of the trivial $G$-module $\Z$.
The cohomology groups $H^*(G,U_d)$ are defined as cohomology groups of the complex
$\Hom_G(\cP,U_d)$ ([Brown, III, \S1]). There exists a chain map of resolutions
$\psi_\cC \colon \cP\to\cC$ (respecting augmentations, i.e., extending by identity on $\Z$
to a chain map of the augmented complexes). Moreover, $\psi_\cC$ is unique up to
chain homotopy ([Brown, I, Lemma 7.4]).
\smallskip
\bf Definition \a.\t.\rm\par
Let $\cC$ be an acyclic $G$-chain complex, $\cP$---a projective resolution of the
trivial $G$-module $\Z$, $\psi_\cC\colon\cP\to\cC$---a chain map of resolutions.
Let $\psi_\cC^*\colon H^d(\Hom_G(\cC,U_d))\to H^d(G,U_d)$ be the map on cohomology induced
by $\psi_\cC$. We define the tautological class:
$$\tau^d_{G,\cC}=\psi^*_\cC(\tau^d_{\cC/G})\in H^d(G,U_d).$$
\smallskip
\rm
These classes are functorial just as the previous ones:
\smallskip
\bf Theorem \a.\t.\sl\par
Let $\cC$ be an acyclic $G$-chain complex, $\cK$---an acyclic $H$-chain complex,
$\phi\colon G\to H$---a group homomorphism, $f\colon\cC\to\cK$---a $\phi$-equivariant chain map.
Consider two maps:
$$
H^d(H,U_d(\cK)))
\mathop{\longrightarrow}\limits^{\phi^*}
H^d(G,U_d(\cK)))
\mathop{\longleftarrow}\limits^{f_*}
H^d(G,U_d(\cC))),
$$
the right one induced by the $f$-induced coefficient map $U_d(\cC)\to U_d(\cK)$.
Then we have:
$$\phi^*\tau^d_{H,\cK}=f_*\tau^d_{G,\cC}.$$
\smallskip
\rm
Proof.
Consider the following diagram:
$$
\matrix{
H^d(\Hom{\!}_H(\cK,U_d(\cK)))
&\mathop{\to}\limits^{f^*}
&H^d(\Hom{\!}_G(\cC,U_d(\cK)))
&\mathop{\leftarrow}\limits^{f_*}
&H^d(\Hom{\!}_G(\cC,U_d(\cC)))\cr
\Big\downarrow
&&
\Big\downarrow
&&
\Big\downarrow\cr
H^d(H,U_d(\cK))
&\mathop{\to}\limits^{\phi^*}
&H^d(G,U_d(\cK))
&\mathop{\leftarrow}\limits^{f_*}
&H^d(G,U_d(\cC))
}
$$
There are tautological classes:
$\tau_{\cK/H}\in H^d(\Hom_H(\cK,U_d(\cK)))$, defined as the class of the
tautological cochain $T_{\cK/G}(k)=[k]$, and a similar
$\tau_{\cC/G}\in H^d(\Hom_G(\cC,U_d(\cC)))$.
Their images in the group
$H^d(\Hom_G(\cC,U_d(\cK)))$ coincide, since both are clearly equal
to the class of $T$ defined by $T(c)=[f(c)]$.
The classes $\tau_{H,\cK}$ and $\tau_{G,\cC}$ are images of $\tau_{\cK/H}$ and $\tau_{\cC/G}$
(respectively) under the vertical maps. Thus, to prove the theorem, we only need to check
that the above diagram is commutative.

Commutativity of the right square:
The vertical maps are induced
by a (unique up to chain homotopy) $G$-map of chain complexes
$P(G)_*\to \cC$. The horizontal maps are induced by the coefficient map $f_*$.
Since these two maps act on different arguments of the ${\rm Hom}$
functor, they commute.

Commutativity of the left square:
That square is the result
of applying a cohomology functor to the following diagram.
$$
\matrix{
\cK
&\mathop{\longleftarrow}\limits^{f_*}
&\cC\cr
\Big\uparrow
&&\Big\uparrow\cr
P(H)_*
&\mathop{\longleftarrow}\limits^{\phi_*}
&P(G)_*
}
$$
The two compositions to compare are $G$-maps from $P(G)_*$ to the acyclic
chain complex $\cK$ (with the $G$-structure induced via $\phi$);
such a map is unique up to chain-homotopy, hence
the compositions are chain-homotopic. After passing to cohomology,
they become equal.\qed
\medskip
\bf Remark \a.\t.\rm\par
The procedure applied in Definition 1.4 to the tautological class works in greater generality, for arbitrary coefficient
groups and arbitrary classes.
In Theorem 4.4 we will need the following version:
Let $\cC$ be an acyclic $G$-chain complex;
$A$---a $G$-module;
$T\in Z^d(\Hom_G(\cC,A))$---an $A$-valued $G$-invariant $d$-cocycle;
$\cP$---a projective resolution of the trivial $G$-module $\Z$;
$\psi_\cC\colon\cP\to\cC$---a chain map of resolutions.
Let $\psi_{\cC}^*\colon \Hom_G(\cC,A)\to \Hom_G(\cP,A)$ be the
cochain map induced by $\psi_{\cC}$.
We define the group cohomology class $\tau\in H^d(G,A)$
associated to $T$ by
$\tau:=[\psi_\cC^*(T)]$.
\medskip

\bf 2. Geometric complexes.\rm
\medskip
\def\a{2}
\m=1
\n=1
Our main source of acyclic $G$-chain complexes is geometry.
Suppose that $G$ acts on an acyclic simplicial complex $X$ by simplicial automorphisms.
Then the simplicial chain complex $C_*X$ is an acyclic $G$-chain complex.
\smallskip
\bf Definition \a.\t.\rm\par
Let $X$ be an acyclic simplicial $G$-complex. The definitions of Section 1 applied to
the 
acyclic simplicial $G$-chain
complex $C_*X$ give rise to:
\item{$\bullet$} the coefficient group $U_d=U_d(X):=U_d(C_*X)$;
\item{$\bullet$} the tautological cocycle $T^d_{X/G}:=T^d_{C_*X/G}$;
\item{$\bullet$} the tautological class $\tau^d_{X/G}:=\tau^d_{C_*X/G}$;
\item{$\bullet$} the tautological group cohomology class $\tau^d_{G,X}:=\tau^d_{G,C_*X}$.
\smallskip\rm
In our considerations the $G$-complexes $X$ will usually arise as
restricted configuration complexes of homogeneous $G$-spaces.
We will typically start from a transitive $G$-action on a space $\P$.
We will use $\P$ as the set of vertices of $X$, and span simplices
of $X$ on tuples of elements of $\P$ satisfying some
genericity conditions. (A typical example: $G=SL(2,K)$, $\P=K^2\setminus\{0\}$,
a tuple of vectors spans a simplex if and only if every two of them are linearly
independent.) This scheme  applies to many algebraic groups
over arbitrary infinite fields.

The acyclicity of these restricted configuration complexes is usually the consequence
of the star-property defined below. 
\smallskip
\bf Definition \a.\t. \rm ([Kr-T])\rm\par
A simplicial complex $X$ has the star-property if for any
finite subcomplex $Y\subseteq X$ there exists a vertex $v\in X^0\setminus Y^0$
joinable with every simplex of $Y$ ($v$ is joinable with a k-simplex $\sigma=[y_0,\ldots,y_k]$
if $v*\sigma=[v,y_0,\ldots,y_k]$ is a $(k+1)$-simplex in $X$). 
\smallskip\rm
\bf Fact \a.\t.\sl\par
If $X$ has the star-property, then it is acyclic.
\smallskip\rm
Proof. Let $z=\sum a_\sigma\sigma$ be  a cycle in $X$. Let $Y$ be
the union of all simplices $\sigma$ that appear in $z$. Let $v$ be a vertex of $X$
witnessing the star-property for $Y$. Then $z=\partial(\sum a_\sigma v*\sigma)$.\qed
\smallskip
For a complex $X$ with the star-property there is another variant of an acyclic chain complex associated to it:
the (non-degenerate) ordered chain complex $C^o_*X$. The group $C^o_kX$ is the free abelian group
whose basis is the set of all $(k+1)$-tuples of vertices of $X$ that span $k$-simplices
(in other words, the set of ordered, non-degenerate $k$-simplices of $X$).
The boundary operator is defined by the usual formula:
$$\partial[v_0,\ldots,v_k]=\sum_{i=0}^k(-1)^i[v_0,\ldots,\widehat{v_i},\ldots,v_k].$$
By the same argument as in Fact \a.3, the complex $C^o_*X$ is acyclic.
(Warning: for finite simplicial complexes the non-degenerate ordered chain complex
does not calculate homology correctly, e.g.~the complex $C^o_*(\Delta^1)$ is not acyclic.)

If a simplicial complex 
is acted upon by a group $G$,
one can use the ordered chain complex to define the coefficient group
$U_d^o:=(C^o_dX/B^o_dX)_G$, the tautological cocycle $T^d_{C_*^oX/G}$ and the tautological class
$\tau^d_{C_*^oX/G}$. (If $X$ has the star-property, one can further define the tautological group
cohomology class $\tau^d_{G,C_*^oX}$.)
There is a natural epimorphic $G$-chain map $C^o_*X\to C_*X$; it induces
an epimorphism $U^o_d\to U_d$. The group $U^o_d$ is usually insignificantly larger
than $U_d$, as we shall see.

The calculations of the groups $U_d$ and $U_d^o$ are often used
in this paper; we now explain how they are done.
Let $X^{(n)}$ be the set of non-degenerate ordered $n$-simplices in a simplicial complex $X$.
Let $R_n$ be a set of representatives of orbits of $G$ on $X^{(n)}$.
For any $\sigma\in X^{(n)}$ we denote by $\sigma_R$ the unique element of $R_n$ that is
$G$-equivalent to $\sigma$. For chains we put $(\sum a_\sigma\sigma)_R=\sum a_\sigma\sigma_R$.
\smallskip
\bf Fact \a.\t.\sl\par
Let $X$ be a simplicial $G$-complex.
\item{a)} The groups $U^o_d$ is the quotient of the free abelian group
with basis $R_d$ by the subgroup spanned by $\{(\partial \rho)_R\mid \rho\in R_{d+1}\}$.
\item{b)} The groups $U_d$ is the quotient of the free abelian group
with basis $R_d$ by the subgroup spanned by
$\{(\partial \rho)_R\mid \rho\in R_{d+1}\}\cup\{(t\rho)_R-\sgn(t)\rho\mid\rho\in R_d,t\in S_{d+1}\}$.
Moreover, in this description one can change the range of $t$ from the permutation group $S_{d+1}$
to any generating set of this group.
\smallskip\rm
The proof is based on the formula $U_d^o=(C_d^oX)_G/\partial(C^o_{d+1})_G$
and an analogous formula for $U_d$. We denote by $c_G$ the image of the chain $c$ in the coinvariants group.
\smallskip
Proof.
We start with a general remark. Suppose that a group $G$ acts on a set $Y$.
Let $\Z[Y]$ be the free abelian group with basis $Y$.
Then $\Z[Y]$ has a natural $G$-module structure, and the coinvariants module
$\Z[Y]_G$ is the free abelian group with basis $Y/G$ (the orbit space of the $G$-action on $Y$).
If $R\subseteq Y$ is a set of representatives of $G$-orbits, then the bijection
$R\ni r\mapsto G\cdot r\in Y/G$ induces the natural isomorphism $\Z[R]\to \Z[Y/G]\to\Z[Y]_G$.

Applying this discussion to the $G$-action on $X^{(n)}$ we see that
$(C_n^oX)_G=\Z[X^{(n)}]_G=\Z[X^{(n)}/G]\simeq\Z[R_n]$. This isomorphism
$(C_n^oX)_G\to\Z[R_n]$ is clearly given by $c_G\mapsto c_R$. Similarly,
$(C_{n+1}^oX)_G$ is isomorphic to $\Z[R_{n+1}]$, which is generated by $R_{n+1}$.
The map $\partial\colon (C_{n+1}^oX)_G\to (C_n^oX)_G$ can be interpreted as
the map $\Z[R_{n+1}]\ni c\mapsto(\partial c)_R\in\Z[R_n]$;
its image is generated by the images of elements of $R_{n+1}$, i.e.~by the set
$\{(\partial\rho)_R\mid \rho\in R_{n+1}\}$. Part (a) is proved.

For part (b): Let $K$ be the kernel of the epimorphism
$C_d^oX\to C_dX$. The group $K$ is generated by
$\{t\sigma-(\sgn{t})\sigma\mid t\in S_{d+1},\sigma\in X^{(d)}\}$
(one can change the range of $t$ from $S_{d+1}$ to any generating set of $S_{d+1}$).
Applying the coinvariants functor to the exact sequence $K\to C_d^oX\to C_dX\to0$
we get the middle row of the following commuting diagram with exact rows and columns. 
$$
\matrix{
&&
(C^o_{d+1}X)_G&
\longrightarrow&
(C_{d+1}X)_G&
\longrightarrow&
0
\cr
&&
\mathop{\Big\downarrow}\limits{\partial}
&&
\mathop{\Big\downarrow}\limits{\partial}
&&
\cr
K_G&\mathop{\longrightarrow}\limits^{\iota_G}&
(C^o_dX)_G&
\longrightarrow&
(C_dX)_G&
\longrightarrow&
0
\cr
&&
\mathop{\Big\downarrow}
&&
\mathop{\Big\downarrow}
&&
\cr
&&
U_d^o
&
\longrightarrow
&
U_d
&&
\cr
&&&&
\Big\downarrow
&&
\cr
&&&&
0
&&
}
$$
A diagram chase shows that an element of $(C^o_dX)_G$ that maps to $0$ in $U_d$
is a sum of images of elements of $K_G$ and $(C^o_{d+1}X)_G$. Consequently,
$$U_d\simeq (C^o_dX)_G/(\partial(C^o_{d+1}X)_G+\iota_GK_G).$$
Therefore, a presentation of
$U_d$ can be obtained from the presentation of $U^o_d$ given in (a) by adjoining extra relations
generating $\iota_GK_G$. These extra relations are images of generators of $K_G$ under $\iota_G$,
i.e.~are of the form $(t\sigma)_G-(\sgn{t})\sigma_G$ ($t\in S_{d+1}$, $\sigma\in X^{(d)}$).
Under the isomorphism $(C^o_dX)_G\to\Z[R_d]$ this form maps to $(t\sigma)_R-(\sgn{t})\sigma_R$.
To finish the proof we will check that $(t\sigma)_R=(t\sigma_R)_R$.
We have $\sigma_R=g\sigma$ for some $g\in G$. This implies that $t\sigma_R=g(t\sigma)$,
and then $(t\sigma_R)_R=(g(t\sigma))_R=(t\sigma)_R$.\qed

\medskip
\bf 3. Boundedness.\rm
\def\a{3}
\m=1
\n=1
\medskip
A group cohomology class in $H^d(G,\R)$ is called bounded if
it can be represented by a bounded cocycle $c\colon S_dBG\to\R$
(or, equivalently, a bounded $G$-invariant $\R$-valued cocycle on $S_dEG$).
Here $S_*BG$ is the singular chain complex of $BG$; a cocycle
$c$ is bounded if there exists $M>0$ such that for each singular
simplex $\sigma\colon\Delta^d\to BG$ we have $|c(\sigma)|\le M$.
Instead of $\R$, one can use other groups with seminorm. In particular,
if $X$ is a simplicial $G$-complex, the coefficient group $U=U_d(X)$
carries a natural seminorm, induced by the $\ell^1$-norm on $C_dX$.
Explicitly, for $u\in U$ we consider all chains
$\sum\alpha_i\sigma_i\in C_dX$ that represent $u$, and we declare
the infimum of $\sum|\alpha_i|$ over all such chains to be $|u|$.
\smallskip
\bf Theorem \a.\t.\sl\par
Suppose that $X$ is an acyclic simplicial $G$-complex with the star-property.
Then the tautological cohomology class $\tau^d_{G,X}\in H^d(G,U)$ is
bounded with respect to the seminorm discussed above.
\smallskip\rm
Proof. We will construct a $G$-chain map $\Psi_*\colon S_*EG\to C_*X$.
For each $n\ge 0$ choose a free basis $\Sigma_n$ of the free $G$-module $S_nEG$.
We define $\Psi_n$ inductively. 
For each $\xi_0\in\Sigma_0$ we choose a vertex
$\Psi_0(\xi_0)\in X^{(0)}$; we extend $\Psi_0$ to $S_0EG$ by $G$-equivariance and linearity.
Once $\Psi_{n-1}$ is defined, we define $\Psi_n$ on $\Sigma_n$ as follows.
For $\xi_n\in\Sigma_n$ we consider $\Psi_{n-1}(\partial\xi_n)=\sum\sigma_i\in C_{n-1}X$.
By the star-property, there exists a vertex $v\in X^{(0)}$ joinable to every
$\sigma_i$; we put $\Psi_n(\xi_n)=\sum v*\sigma_i$, so as to have
$\partial\Psi_n(\xi_n)=\Psi_{n-1}(\partial\xi_n)$. Then we extend $\Psi_n$
to $S_nEG$ by $G$-equivariance and linearity. A straightforward induction shows
that for any singular simplex $\xi_n\in S_nEG$ the chain $\Psi_n(\xi_n)$
is a sum of at most $(n+1)!$ simplices.

The class $\tau^d_{G,X}$ is represented by the cocycle $T^d_{X/G}\circ\Psi_d$.
The tautological cocycle $T^d_{X/G}$ has norm at most 1---it maps
a simplex to its class in $U_d$, and that class has norm $\le1$
by definition of the seminorm. Therefore, for any singular simplex $\sigma_d$ in $EG$ we have
$$|T^d_{X/G}(\Psi_d(\sigma_d))|\le (d+1)!.$$
\qed

\bf Remark \a.\t.\rm\par
There is a different approach to bounded group cohomology, based on the
standard homogeneous resolution of the trivial $G$-module $\Z$
(cf.~[Brown, I, \S5]).
That approach is equivalent to the one used above, as shown in [Gro, pp.~48--49];
for a more detailed account see [L\"oh, 2.5.5]. In these references real coefficients
are used, but the proof works for coefficients in an arbitrary abelian group with seminorm.
\medskip
\bf 4. Characteristic classes.\rm
\def\a{4}
\m=1
\n=1
\medskip
A cohomology class $\alpha$ of a (discrete) group $G$ can serve as a
characteristic class of (flat) $G$-bundles. Suppose that $\alpha$ is obtained
from a $G$-invariant cocycle on an acyclic $G$-space $X$ as in Remark 1.6. 
Then it is possible to describe the characteristic class using the cocycle directly, by-passing $\alpha$ (see Theorem \a.4). 
This section is organized as follows. We start by recalling the connection between group
cohomology and characteristic classes. Next, we 
describe the classical de Rham version of characteristic classes of flat bundles.
Then we discuss auxiliary notions and notation and, finally, we state
and prove the main statement, Theorem \a.4. 
(Recall that
we consider $G$ with discrete topology, so that all $G$-bundles are flat---with locally constant
transition functions---and $BG$ is $K(G,1)$.)\rm

Let $\alpha\in H^d(G,A)=H^d(BG,A)$ be a cohomology class of a group $G$.
The space $BG$ is the base of a universal principal $G$-bundle $EG$.
Every principal $G$-bundle $P$ over a (paracompact) base space $B$
has a classifying map:
a map $f_P\colon B\to BG$ such that $f_P^*EG\simeq P$. 
The map $f_P$ is unique up to homotopy. Notice that we use $f^*\xi$
to denote the pull-back of the bundle $\xi$ via the map $f$, and we also
use $f^*\tau$, $f^*T$  for the pull-back of a cohomology class $\tau$
or of a cocycle $T$. Though occasionally confusing, this dual usage is
standard practice in bundle theory.

\smallskip
\bf Definition \a.\t.\rm\par
The cohomology class
$\alpha(P):=f_P^*(\alpha)\in H^d(B,A)$
is functorial in $P$,
and is called the characteristic class (corresponding to $\alpha$)
of the bundle $P$.
\smallskip\rm
\def\base{B}
In this definition the $G$-module $A$ may have non-trivial $G$-structure.
Then the groups $H^d(G,A)$ and $H^d(\base,A)$ are cohomology groups with twisted coefficients,
i.e.~with coefficients in a flat $G$-bundle (local system) with fibre $A$.
For $H^d(G,A)=H^d(BG,A)$ the bundle is $EG\times_GA$; for $H^d(\base,A)$ we use
$P\times_GA$.
We have
$P\times_GA=f_P^*(EG\times_GA)$, so that the coefficient system
used over $BG$ pulls back to the one used over $\base$; therefore we get
a map $f_P^*\colon H^d(G,A)\to H^d(\base,A)$.
\smallskip\rm
In 
de Rham theory there is a construction of characteristic classes of flat bundles
that does not explicitly refer to $BG$. In fact, it gives an explicit cocycle representative
of the characteristic cohomology class in terms of a section.
Suppose that $P\to \base$ is a principal flat $G$-bundle over a manifold $\base$, and that
$\omega\in \Omega^d(X)$ is a $G$-invariant closed form on a contractible $G$-manifold $X$.
To these data we will associate a class in $H^d_{\rm DR}(\base)$. 
We start by forming the associated bundle $E=P\times_GX$ with fibre $X$.
Then we choose a section $s\colon \base\to E$; it exists and is homotopically unique
because $X$ is contractible.
Now the idea is that a section $s$ of a flat bundle is an
ill-defined---$G$-ambivalent---map
from the base to the fibre. The $G$-ambivalence is countered by the $G$-invariance of $\omega$,
so that the pull-back of $\omega$ by $s$ is well-defined.
Let us be more precise.
Let $\varphi_U\colon E|_U\to U\times X$ be local trivializations of $E$.
Composing $\varphi_U$ with ${\rm pr}_2\colon U\times X\to X$ we get a map
$\psi_U\colon E|_U\to X$. The compositions $\psi_U\circ s|_U\colon U\to X$
are locally defined maps; these maps are not compatible.
However, due to the $G$-invariance of $\omega$,
the forms $\omega_U=(\psi_U\circ s|_U)^*\omega\in \Omega^d(U)$ are compatible
and define a global closed form in $\Omega^d(\base)$. Slightly abusing the notation
we denote this form by $s^*\omega$.
The cohomology class of $s^*\omega$ in $H^d_{\rm DR}(\base)$
is a characteristic class of the bundle $P$.
An alternative description
is to define the global form $\omega^E$ on $E$ by gluing the compatible collection
of forms $\psi_U^*\omega\in \Omega^d(E|_U)$, and then take $s^*\omega^E$ in the standard sense.  
(See [Morita, Chapter 2] for more information on these classes.)

Let us pass to the simplicial setting. Let $P\to\base$ be a principal
$G$-bundle over a $\Delta$-complex $\base$. (For a basic discussion of $\Delta$-complexes see [Hatcher, Section 2.1].)
Let $T\in Z^d(\Hom_G(C_*X,A))$ be an $A$-valued $G$-invariant simplicial
cocycle on an acyclic simplicial $G$-complex $X$, and let $\tau\in H^d(G,A)$ be the associated
cohomology class (as in Remark 1.6).
The characteristic class of $P$ (corresponding to $\tau$) is the cohomology class
$\tau(P)\in H^d(\base,A)$ (see Definition \a.1).
We will use the strategy explained in the de Rham setting
and obtain a cochain on $\base$ representing $\tau(P)$ (see Theorem \a.4).

To deal with sections in the simplicial context we introduce a special family of trivializations.
Let $P\to \base$ be a principal $G$-bundle over a $\Delta$-complex $\base$.
Let $X$ be a simplicial $G$-complex. Let $E=P\times_GX$ be the associated bundle
with fibre $X$. Consider a simplex $\sigma\colon\Delta\to\base$, part of the $\Delta$-complex structure.
The bundle $\sigma^*P$ is a flat principal $G$-bundle over a simplex, hence
it has flat sections. Any such flat section $r\colon\Delta\to \sigma^*P$
induces a trivialization of $\sigma^*E\simeq\sigma^*P\times_GX$---the map
$$\Delta\times X\ni(p,x)\mapsto [r(p),x]\in \sigma^*P\times_GX$$
is an isomorphism, whose inverse $\varphi_{\sigma,r}$ is a trivialization.
We put $\psi_{\sigma,r}={\rm pr}_2\circ\varphi_{\sigma,r}\colon \sigma^*E\to X$.
Notice that all possible flat sections of $\sigma^*P$ are $G$-related, and that
$$\psi_{\sigma,rg}=g^{-1}\psi_{\sigma,r}.\leqno(\a.\f)$$
Moreover, if $\sigma_i$ is a face of $\sigma$
(say $\sigma_i=\sigma|_{\Delta(i)}$, where $\Delta(i)={[e_0,\ldots,\widehat{e_i},\ldots,e_n]}$), then
$$\psi_{\sigma,r}|_{\sigma_i^*E}=\psi_{\sigma_i,r|_{\Delta(i)}}.\leqno(\a.\f)$$
We will now use the maps $\psi_{\sigma,r}$ to define simplicial sections.
\smallskip
\bf Definition \a.\t.\rm\par
Let $\base$ be a $\Delta$-complex, $X$---a simplicial $G$-complex, $P\to\base$---a principal $G$-bundle,
$E=P\times_GX$---the associated bundle over $\base$ with fibre $X$.
A section $s\colon\base\to E$ is called simplicial, if for every simplex
$\sigma\colon \Delta\to \base$ from the $\Delta$-structure of $\base$, and for any
$\psi_{\sigma,r}\colon \sigma^*E\to X$ as described above, the composition
$\psi_{\sigma,r}\circ s\circ \sigma\colon \Delta\to X$ is an affine map of $\Delta$
onto some simplex of the simplicial structure of $X$---possibly
onto a simplex of dimension smaller than $\dim\Delta$
(the composition $s\circ\sigma$ defines a section of $\sigma^*E$, because for $p\in\Delta$
we have $(\sigma^*E)_p=E_{\sigma(p)}$).
\smallskip\rm
\bf Remark.\rm\par
A simplicial section in uniquely determined by its values at the vertices
of the base.
\smallskip\rm
A twisted cochain in $C^d(\base,A)$ assigns to a simplex $\sigma\colon\Delta\to \base$
a value in $(P\times_GA)_{\sigma(e_0)}$---the fibre of the coefficient bundle
over the initial vertex of $\sigma$. This value extends to a (unique) flat section
of $\sigma^*(P\times_GA)=\sigma^*P\times_GA$. A flat section of that bundle
can be described as $[r,a]$, where $r$ is a section of $\sigma^*P$ and $a\in A$;
for each $g\in G$ the  pair $[rg,g^{-1}a]$ defines the same section, therefore one can also
describe sections as continuous (locally constant) $G$-maps $\sigma^*P\to A$---or $G$-maps
from the $G$-torsor of flat sections of $\sigma^*P$ to $A$.
\smallskip
\bf Definition \a.\t.\rm\par
Let $\base$ be a $\Delta$-complex, $X$---a simplicial $G$-complex, $P\to\base$---a principal $G$-bundle,
$E=P\times_GX$---the associated bundle over $\base$ with fibre $X$,
$s$---a simplicial section of $E$.
Consider a simplex $\sigma\colon\Delta\to\base$ from the $\Delta$-structure of $\base$
and flat sections $r$ of $\sigma^*P$.
Then the expression
$T(\psi_{\sigma,r}\circ s\circ\sigma)$ is $G$-equivariant in $r$
(due to (\a.1) and the fact that $T$ is a $G$-map). The formula
$$s^*T(\sigma)=[r,T(\psi_{\sigma,r}\circ s\circ\sigma)]$$
defines the cochain $s^*T\in C^d(\base,A)$ (with twisted coefficients).
\smallskip\rm
(The image of the map
$\psi_{\sigma,r}\circ s\circ\sigma$ is a simplex in $X$, on which we put the
orientation corresponding under this map to the standard orientation of the standard simplex;
we interpret the argument of $T$ as that oriented simplex.
If the image of $\psi_{\sigma,r}\circ s\circ\sigma$ has dimension smaller than $d$,
we interpret the argument of $T$ as the zero chain.)
\smallskip\rm
\bf Remark.\rm\par
The fact that $s^*T$ is a cocycle will follow from the proof of the next theorem. 
\smallskip
\bf Theorem \a.\t.\sl\par
Let $T$ be an $A$-valued $G$-invariant cocycle on an acyclic simplicial $G$-complex $X$.
Let $\tau\in H^d(G,A)$ be the associated group cohomology class (as in Remark 1.6).
Let $P\to\base$ be a principal (flat) $G$-bundle over a $\Delta$-complex $\base$.
Let $s\colon \base\to P\times_GX$ be a simplicial section. Then the class
$\tau(P)\in H^d(\base,A)$---the characteristic class of $P$ corresponding to $\tau$---is represented
by the simplicial cocycle $s^*T\in Z^d(\base,A)$.
\smallskip\rm
Proof. 

The total space $EG$ of the universal principal $G$-bundle $EG\to BG$ is contractible
(it is also the universal cover of $BG$).
The $G$-action on $EG$ is free.
Therefore, the singular chain complex $S_*EG$ is a projective (in fact, free) resolution
of the trivial $G$-module $\Z$. Moreover, we have $(S_*EG)_G\simeq S_*BG$.
Let $\Psi=\Psi_{C_*X}\colon S_*EG\to C_*X$ be a resolution map from $S_*EG$ 
to the simplicial chain complex of $X$.
This map induces the map $\Psi^*\colon \Hom_G(C_*X,A)\to \Hom_G(S_*EG,A)$, and
$\tau=[\Psi^*T]\in H^d(\Hom_G(S_*EG,A))=H^d(G,A)$.

Let $f=f_P\colon\base\to BG$ the a classifying map of the bundle $P$,
and let $F\colon P\to EG$ be a $G$-bundle map covering $f$.
Then $\tau(P)=f^*\tau=[f^*\Psi^*T]$. Let us describe the cocycle $f^*\Psi^*T$ explicitly.
This cocycle should assign to any simplex $\sigma\colon\Delta\to \base$ (from the $\Delta$-structure of $\base$)
a value in $(P\times_GA)_{\sigma(e_0)}$; as explained in the paragraph
preceding Definition \a.3, the choice of that value
is equivalent to the choice of a $G$-map from the $G$-set of flat sections $r$ of $\sigma^*P$ to the $G$-module $A$.
Suppose 
that $r\colon\Delta\to \sigma^*P$ is a (flat) section.
Then $F\circ r$ is a singular simplex in $EG$.
The map $r\mapsto T(\Psi(F\circ r))\in A$
is a $G$-map (since each of $F$, $\Psi$, $T$ is a $G$-map);
it defines 
the value of the cochain $f^*\Psi^*T$ (representing $\tau(P)$) on the simplex $\sigma$:
$$f^*\Psi^*T(\sigma)=[r,T(\Psi(F\circ r))].\leqno(\a.\f)$$

The cochain $f^*\Psi^*T$ depends on several choices:
\item{(1)} one can choose the space $BG$---within the homotopy type;
\item{(2)} one can choose $f\colon\base\to BG$---within the homotopy class;
\item{(3)} one can choose $\Psi$---all resolution maps are possible.

Our strategy is to exploit these choices to ensure that $f^*\Psi^*T=s^*T$.

Let $\Sigma(\base)$ be the set of all simplices $\sigma$ forming the $\Delta$-structure of $\base$.
\smallskip
\bf Lemma \a.\t.\par\sl
One can choose the space $BG$ and the map $f$ so that:
\item{(a)} $f^*EG\simeq P$ (i.e.~$f$ is a classifying map of $P$);
\item{(b)} all the maps $f\circ\sigma$ for $\sigma\in\Sigma(\base)$ are pairwise
distinct.
\smallskip\rm
Proof. For dimension $d$ let $m_d$ be the barycentre of the standard simplex $\Delta^d$.
For each $\sigma\in\Sigma(\base)$ we put
$p_\sigma=\sigma(m_{\dim{\sigma}})$. Then we choose a collection
of pairwise different points $(x_\sigma)_{\sigma\in\Sigma(\base)}$ in $BG$.
(If $BG$ is too small for that, we change it by wedging it with a contractible
space of sufficiently large cardinality.)
Finally, we perform a homotopy of $f$ (inductively over skeleta) to ensure $f(p_\sigma)=x_\sigma$.
\qed(Lemma)
\smallskip
Let $f\colon\base\to BG$ be a classifying map of $P$ satisfying the conditions of Lemma \a.5.
Let $F\colon P\to EG$ be a $G$-bundle map covering $f$ (the composition of an
isomorphism $P\to f^*EG$ with the canonical map $f^*EG\to EG$).
\smallskip
\bf Lemma \a.\t.\sl\par\nobreak
One can choose the resolution map $\Psi\colon S_*EG\to C_*X$ so that $f^*\Psi^*T=s^*T$.
\smallskip\rm
Proof. Let us discuss how $\Psi$ may be constructed. 
For each $n\ge 0$ choose a free
basis $\Sigma_n$ of the free $G$-module $S_nEG$. Define $\Psi_n$ inductively.
The base case is $\Psi_{-1}={\rm Id}_\Z\colon\Z\to\Z$, with $\Z$ connected
to the resolutions by the augmentation maps $\partial\colon S_0EG\to\Z$,
$\partial\colon C_0X\to\Z$.
Once
$\Psi_{n-1}$ is defined, calculate---for every $\sigma\in\Sigma_n$---the cycle $\Psi_{n-1}(\partial\sigma)$.
Since $C_*X$ is acyclic, this cycle is a boundary of some $n$-chain; pick one such chain
and define it to be $\Psi_n(\sigma)$.
A crucial remark is that if, for some $\eta\in\Sigma_n$ and some $n$-simplex $\xi$ in $X$, we have
$\Psi_{n-1}(\partial\eta)=\partial\xi$, then we may put $\Psi_n(\eta)=\xi$.
Once $\Psi_n$ is defined on $\Sigma_n$, we extend it to $S_nEG$ by $G$-equivariance and linearity.

For each $(\sigma\colon\Delta\to\base)\in\Sigma(\base)$ choose a flat section $r(\sigma)\colon\Delta\to\sigma^*P$.
Composing this section with the canonical bundle map $\sigma^*P\to P$, and then with $F\colon P\to EG$, we get
a singular simplex $F\circ r(\sigma)$ in $EG$. We denote this simplex by
$\widetilde{f\circ\sigma}$---it is a lift of $f\circ\sigma$.
All the lifts $\widetilde{f\circ\sigma}$ are pairwise $G$-inequivalent, because all $f\circ\sigma$ are pairwise distinct.
Therefore we may choose the free bases $\Sigma_n$ so that they contain all the lifts $\widetilde{f\circ\sigma}$
(for $\sigma\in\Sigma(\base)$). We would like to define:
$$\Psi(\widetilde{f\circ\sigma})=\psi_{\sigma,r(\sigma)}\circ s\circ\sigma.\leqno(\a.\f)$$
To be able to do that we need to check that $$\Psi(\partial(\widetilde{f\circ\sigma}))=
\partial(\psi_{\sigma,r(\sigma)}\circ s\circ\sigma).\leqno(\a.\f)$$
Let $\sigma_i=\sigma|_{\Delta(i)}$,
where $\Delta(i)={[e_0,\ldots,\widehat{e_i},\ldots,e_n]}$, $n=\dim{\sigma}$. (Strictly speaking, we should
also use an extra map identifying $\Delta(i)$ with the standard simplex. We will ignore this in order not to overburden
the notation.)
We have:
$$\partial(\widetilde{f\circ\sigma})=\sum_{i=0}^n(-1)^i(\widetilde{f\circ\sigma})|_{\Delta(i)}.
\leqno(\a.\f)$$
Observe that $(\widetilde{f\circ\sigma})|_{\Delta(i)}$ is a lift of $f\circ(\sigma|_{\Delta(i)})$;
therefore $(\widetilde{f\circ\sigma})|_{\Delta(i)}
=(\widetilde{f\circ\sigma_i})\cdot g(i)$ for some
$g(i)\in G$. By induction on the dimension we know that
$$ \Psi((\widetilde{f\circ\sigma_i})\cdot g(i))=g(i)^{-1}\Psi(\widetilde{f\circ\sigma_i})
=g(i)^{-1}(\psi_{\sigma_i,r(\sigma_i)}\circ s\circ\sigma_i)
=\psi_{\sigma_i,r(\sigma_i)g(i)}\circ s\circ\sigma_i,
$$
the last equality following from (\a.1). We may finally write
$$\Psi(\partial(\widetilde{f\circ\sigma}))=\sum_{i=0}^n (-1)^i\psi_{\sigma_i,r(\sigma_i)g(i)}\circ s\circ\sigma_i.\leqno(\a.\f)$$
On the other hand,
$$\partial(\psi_{\sigma,r(\sigma)}\circ s\circ\sigma)=\sum_{i=0}^n(-1)^i\psi_{\sigma,r(\sigma)}\circ s\circ\sigma_i.\leqno(\a.\f)$$
Notice that $s\circ\sigma_i$ is a section of $\sigma_i^*E$ (where $E=P\times_GX$); therefore (\a.2) applies and yields
$$\psi_{\sigma,r(\sigma)}\circ s\circ\sigma_i=\psi_{\sigma_i,r(\sigma)|_{\Delta(i)}}\circ s\circ\sigma_i.\leqno(\a.\f)$$
The flat sections $r(\sigma)|_{\Delta(i)}$ and $r(\sigma_i)$ of the $G$-bundle $\sigma_i^*P$
are $G$-related: $r(\sigma)|_{\Delta(i)}=r(\sigma_i)\cdot g$ for some $g\in G$. Recall that,
for any $\eta\in\Sigma(\base)$, the singular simplex $\widetilde{f\circ\eta}$ is the composition
of $r(\eta)$ with a $G$-bundle map $\eta^*P\to EG$. It follows that
$(\widetilde{f\circ\sigma})|_{\Delta(i)}=(\widetilde{f\circ\sigma_i})\cdot g$;
therefore $g=g(i)$. Consequently,
$$\psi_{\sigma_i,r(\sigma)|_{\Delta(i)}}\circ s\circ\sigma_i=\psi_{\sigma_i,r(\sigma_i)g(i)}\circ s\circ\sigma_i.\leqno(\a.\f)$$
Putting (\a.8), (\a.9) and (\a.10) together we get
$$\partial(\psi_{\sigma,r(\sigma)}\circ s\circ\sigma)=\sum_{i=0}^n(-1)^i\psi_{\sigma_i,r(\sigma_i)g(i)}\circ s\circ\sigma_i.\leqno(\a.\f)$$
Comparing this with (\a.7) we obtain (\a.5). We may therefore define $\Psi$ so that (\a.4) holds for all $\sigma\in\Sigma(\base)$.
Then, for $d$-dimensional $\sigma$ we get
$$\eqalign{
f^*\Psi^*T(\sigma)
&=[r(\sigma),T(\Psi(F\circ r(\sigma)))]\cr
&=[r(\sigma),T(\Psi(\widetilde{f\circ\sigma}))]
=[r(\sigma),T(\psi_{\sigma,r(\sigma)}\circ s\circ\sigma)]=s^*T(\sigma).
}\leqno(\a.\f)$$
\qed(Lemma \a.6)

Lemma \a.6 implies the theorem: $\tau(P)$ is represented by $f^*\Psi^*T$, which is equal to $s^*T$.

\qed(Theorem \a.4)
\smallskip
\bf Remark.\rm\par
In our applications of Theorem \a.4 the coefficients will be either untwisted or
only mildly twisted (e.g.~a $GL(2,K)$-module $A$ which is trivial as an $SL(2,K)$-module).

\medskip
\bf 5. Homological core.\rm
\medskip
\def\a{5}
\m=1
\n=1
Consider an acyclic simplicial $G$-complex $X$, the associated coefficient group $U_d$
and the tautological cohomology classes
$\tau_{X/G}\in H^d(\Hom_G(C_*X,U_d))$ and $\tau_{G,X}\in H^d(G,U_d)$.
\smallskip
\bf Question: \rm Is it possible to represent these tautological classes
by cocycles with coefficients in a proper subgroup of $U_d$?
\smallskip
In general, there is a candidate subgroup.
The coefficient group $U_d=(C_dX)_G/\partial(C_{d+1}X)_G$ has a natural homomorphism
$\partial\colon U_d\to(C_{d-1}X)_G$ (induced by the usual $\partial\colon C_dX\to C_{d-1}X$).
\smallskip
\bf Definition \a.\t.\par\rm
The homological core $hU_d$ of the group $U_d$ is the kernel of the map $\partial\colon U_d\to (C_{d-1}X)_G$.
\smallskip
The following theorem states a weaker property then asked for above, but is quite general.
\rm
\smallskip
\bf Theorem \a.\t.\sl\par
Let $X$ be an acyclic simplicial $G$-complex. Let  $\tau_{G,X}\in H^d(G,U_d)$
be the associated tautological class, and let $z\in H_d(BG,\Z)$ be a homology class.
Then $\langle\tau_{G,X},z\rangle\in hU_d$.
\smallskip
\rm
Proof.
The map 
$\partial\colon U_d\to (C_{d-1}X)_G$
of coefficient groups 
induces horizontal maps in the commutative diagram:
$$
\matrix{
H^d(\Hom_G(C_*X,U_d))
&
\mathop{\to}\limits^{\partial_*}
&
H^d(\Hom_G(C_*X,(C_{d-1}X)_G))
\cr
\Big\downarrow
&&
\Big\downarrow
\cr
H^d(G,U_d)
&
\mathop{\to}\limits^{\partial_*}
&
H^d(G,(C_{d-1}X)_G).
}
$$

The class $\tau_{X/G}\in H^d(\Hom_G(C_*X,U_d))$ is mapped to $0$ by $\partial_*$.
Indeed, the class $\tau_{X/G}$ is represented by
the tautological cocycle given by $T_{X/G}(\sigma^d)=[\sigma^d]$.
Let $t\in \Hom_G(C_{d-1},(C_{d-1})_G)$ also be tautological:
$t(\sigma^{d-1})=[\sigma^{d-1}]$. Then
$$(\partial_*T_{X/G})(\sigma^d)=\partial[\sigma^d]=[\partial\sigma^d]=t(\partial\sigma^d)=(\delta t)(\sigma^d),$$
that is, $\partial_*T_{X/G}=\delta t$, hence $\partial_*\tau_{X/G}=0$.

It now follows from the diagram that $\partial_*\tau_{G,X}=0$ as well.
Therefore, $\partial\langle\tau_{G,X},z\rangle=\langle\partial_*\tau_{G,X},z\rangle=0$.\qed

\bigskip
\vfill\eject
\centerline{\bf II. $GL(2,K)$}
\bigskip
In this part we describe some results of Nekov\'a$\check{\rm r}$ (cf.~[Ne]) from our point of view.
This provides a perfect illustration of the general method.
\medskip
\bf 6. Review of possible actions.\rm
\medskip
\def\a{6}
\m=1
\n=1
The first example where the general approach from part I gives something interesting
is when $G=GL(2,K)$. To proceed we need a simplicial action of $G$ on a complex $X$
with desired properties, one of them being high transitivity.
Such $X$ can be constructed by taking as the vertex set a homogeneous space $G/S$
for some $S$ and studying the notion of ``generic $k$-tuple''. In $GL(2,K)$ we have the
following interesting subgroups.
\item{1)} {$S=\pmatrix{1&*\cr 0&*}$. Then $G/S=K^2\setminus\{0\}$. A tuple
is generic if it consists of pairwise linearly independent vectors.
The action of $GL(2,K)$ is effective and transitive on generic pairs. We will discuss this case
later.}
\item{2)} {$S=\pmatrix{*&0\cr 0&*}$. Then $G/S=(\P^1(K)\times\P^1(K))\setminus\Delta$ (where
$\Delta$ is the diagonal). Generic $k$-tuples  are tuples of pairs $(p_i,q_i)$ with all
the points $p_i$ and $q_j$ distinct. Here even the action on pairs is not transitive, because
the cross-ratio $(p_1,q_1,p_2,q_2)$ is preserved. The action of $GL(2,K)$ is not effective,
it factors through $PGL(2,K)$.}
\item{3)} {$S=\pmatrix{*&*\cr 0&*}$. Then $G/S=\P^1(K)$. The action factors
through $PGL(2,K)$. Generic $k$-tuples are tuples of distinct projective points.
The action is triply transitive.}
\medskip
In the first two examples our approach yields very big groups of coefficients.
We can compute them (and we do so in the $S=\pmatrix{1&*\cr0&*}$ case), but
we cannot say much about them. In the third case transitivity is higher,
hence the coefficient group is smaller and easier to understand. We do the computation
in detail (following Nekov\'a$\check{\rm r}$).
Actually in this case something interesting happens: while $PGL(2,K)$ acts transitively
on triples, the large normal subgroups $PSL(2,K)$ acts transitively on pairs, while
its orbits on triples are indexed by $\sqcl$---the group of square classes.
This gives an untwisted cohomology class for $PSL(2,K)$ and a (slightly) twisted
cohomology class for $PGL(2,K)$.

Now we proceed with the description of the $S=\pmatrix{1&*\cr 0&*}$ case. We do not discuss
the $S=\pmatrix{*&0\cr0&*}$ case.
\smallskip
\bf $GL(2,K)$-action on non-zero vectors. \rm

Consider the $GL(2,K)$-action on $GL(2,K)/\pmatrix{1&*\cr0&*}\simeq K^2\setminus\{0\}$.
We declare a tuple of non-zero vectors in $K^2$ to be generic if its elements are
pairwise linearly independent. The complex $X$ has $k$-simplices spanned on generic $(k+1)$-tuples.

The action of $GL(2,K)$ is transitive on $1$-simplices. However, the action of its
large normal subgroup $SL(2,K)$ does have an invariant: $(v,w)\mapsto\det(v,w)\in\dotK$.
Thus, we have a potentially nontrivial $1$-dimensional cocycle, with coefficients in a quotient
of $\Z[\dotK]$. However, any given generic
triple can be normalized by an element of $SL(2,K)$ to
$$\left({1\choose 0}, {0\choose a},{-x/a\choose y}\right)\leqno(\a.\f)$$
with non-zero $a$, $x$ and $y$;
the three determinants are $a$, $y$, $x$. To get the coefficient group
we have to divide $\Z[\dotK]$ by relations $\symb(a)-\symb(x)+\symb(y)$, for
all non-zero  $a$, $x$ and $y$. The resulting coefficient group is trivial,
hence so is the cocycle.

The action of $GL(2,K)$ on generic triples of vectors
has a complete invariant in $\dotK\times\dotK\times\dotK$:
the triple of pairwise determinants. The resulting coefficient group
is the quotient of $\Z[\dotK\times\dotK\times\dotK]$ by the following
relation:
$$({be-cd\over a},e,c)-({be-cd\over a},d,b)+(e,d,a)-(c,b,a)=0$$
(all entries assumed non-zero). We skip the details, as they are not dissimilar
to ones in the calculation presented later, and the result is not especially meaningful.
\smallskip
\bf $GL(2,K)$-action on projective line. \rm

Consider the $GL(2,K)$-action on $GL(2,K)/\pmatrix{*&*\cr0&*}\simeq \P^1(K)$.
We declare a $(k+1)$-tuple of projective points generic if they are pairwise distinct;
we span $k$-simplices on such tuples. Thus, the complex $X$ is the (infinite) simplex
with vertex set $\P^1(K)$.

The $GL(2,K)$ (in fact, $PGL(2,K)$) action on $X$ is transitive on $2$-simplices.
However, the $SL(2,K)$-action on $2$-simplices has an invariant with values in $\sqcl$---the set
of square classes.
Our procedure will produce a cocycle with constant coefficients
(in a quotient group of $\Z[\sqcl]$) for $PSL(2,K)$, and with twisted coefficients
for $PGL(2,K)$. We discuss this in detail in the next section.

The cross-ration is a complete invariant of (ordered) $3$-simplices, i.e~of
$4$-tuples of distinct points in $\P^1(K)$, under the action of
$PGL(2,K)$.
Thus, our approach yields a $3$-cocycle
with coefficients in $\Z[\dotK\setminus\{1\}]/I$, where $I$ is the subgroup
spanned by
$$\left[{\lambda(\mu-1)\over\mu(\lambda-1)}\right]-
\left[{\mu-1\over \lambda-1}\right]+
\left[{\mu\over\lambda}\right]-
[\mu]+[\lambda]
$$
with $\lambda,\mu\in \dotK\setminus\{1\}$, $\lambda\ne\mu$.
This is related to the dilogarithm function; we do not pursue it further
(but see, e.g., [BFG]).
\medskip
\bf 7. Action on triples of points in the projective line.\rm
\medskip
\def\a{7}
\m=1
\n=1
In this section we consider the action of $G=PSL(2,K)$ on the projective line
$\P=\P^1(K)$, for an infinite field $K$. We define $X$ as the (infinite) simplex
with vertex set $\P$; in other words, we span simplices of $X$ on tuples of generic
(i.e.~pairwise distinct) points in $\P$. The complex $X$ has the star-property
and is contractible. To the induced action of $G$ on $X$ the formalism of Part I
(see Definition 2.1) associates the coefficient group $U_2$, the tautological cocycle $T$,
and a tautological cohomology class of $G$.
\smallskip
\bf Theorem \a.\t.\sl\par
The group $U_2$ is isomorphic to $W(K)$, the Witt group of the field $K$.
\rm
\smallskip
Let us briefly recall the definition of $W(K)$
(for details see [EKM, Chapter I]). The isometry classes of non-degenerate
symmetric bilinear forms over $K$ form a semi-ring, with direct sum as addition
and tensor product as multiplication. Passing to the Grothendieck group
of the additive structure of this semi-ring we obtain the
Witt--Grothendieck ring $\widehat{W}(K)$. The Witt ring $W(K)$ is the quotient
of $\widehat{W}(K)$ by the ideal generated by the hyperbolic plane---the
form with matrix $\left({0\,1\atop 1\,0}\right)$.
Both rings have explicit presentations in terms of generators and relations,
cf.~[EKM, Theorems 4.7, 4.8].
\smallskip
Proof of Theorem \a.1. We apply Fact 2.4.

\bf Generators. \rm  We need to find the orbits of the $G$-action
on the set of generic triples of points in $\P$. We denote by $[v]$ the point in
$\P$ determined by the vector $v\in K^2$; for $v={a\choose b}$ we shorten
$\left[{a\choose b}\right]$ to ${a\broose b}$.

\smallskip
\bf Lemma \a.\t.\sl\par
Every generic triple $([u],[v],[w])$ of points in $\P$ is $G$-equivalent to
a triple of the form
$$t_\lambda=\left({1\broose0},{0\broose1},{1\broose\lambda}\right),
\qquad \hbox{\rm where}\quad \lambda={\det(u,v)\det(v,w)\det(w,u)}.$$
Triples $t_\lambda$, $t_\mu$ are equivalent if and only if $\lambda/\mu\in\dotK^2$ (the set of squares in $\dotK$). 
\rm\par
Proof. There exists $g\in SL(2,K)$ such that $gu={1\choose0}$, $gv={0\choose\gamma}$ for $\gamma=\det(u,v)$.
Then $gw={\alpha\choose\beta}=\alpha{1\choose\beta/\alpha}$ for some $\alpha,\beta\in\dotK$, so that
$$g\colon([u],[v],[w])\mapsto
\left({1\broose0},{0\broose\gamma},{\alpha\broose\beta}\right)=
\left({1\broose0},{0\broose1},{1\broose\lambda}\right)$$
for $\lambda=\beta/\alpha$. Notice that
$$\eqalign{
\lambda={\beta\over\alpha}=
\left|\matrix{1&0\cr0&\gamma
}\right|\cdot
\left|\matrix{0&\alpha\cr\gamma&\beta
}\right|^{-1}\cdot
\left|\matrix{\alpha&1\cr\beta&0
}\right|&=
{\det(gu,gv)\det(gv,gw)^{-1}\det(gw,gu)}
\cr
&=
{\det(u,v)\det(v,w)^{-1}\det(w,u)}.}\leqno(\a.\f)$$
Notice that the stabilizer of
$\left({1\broose0},{0\broose1}\right)$
in $SL(2,K)$ consists of diagonal matrices of the form $\pmatrix{\alpha^{-1}&0\cr0&\alpha}$.
Such matrix maps $t_\lambda$ to
$$\left({\alpha^{-1}\broose0},
{0\broose\alpha},
{\alpha^{-1}\broose\alpha\lambda}
\right)=
\left({1\broose0},
{0\broose1},
{1\broose\alpha^2\lambda}
\right)
=t_{\alpha^2\lambda}.$$
The last claim of the lemma follows.
Finally, the class in $\sqcl$ of $\lambda$ given by (\a.1)
is the same as the class of
$\det(u,v)\det(v,w)\det(w,u)$.
\qed{(Lemma)}
\smallskip
The Lemma and Fact 2.4 imply that $U_2$ is the quotient of $\Z[\sqcl]$ by two sets of relations
(boundary relations and alternation relations). The generator of $\Z[\sqcl]$ corresponding to
$t_\lambda$ (and the image of this generator in $U_2$) will be denoted by $[\lambda]$
and called the symbol of the triple.
\smallskip
\bf Alternation relations. \rm
The transposition of the first two points of a triple maps $t_\lambda$ to
$$
\left(
{0\broose1},
{1\broose0},
{1\broose\lambda}
\right),$$
which can be transformed by $\pmatrix{0&1\cr-1&0}$ to
$$\left({1\broose0},{0\broose-1},{\lambda\broose-1}\right)=
\left({1\broose0},{0\broose1},{1\broose-\lambda^{-1}}\right)
.$$
Since $[-\lambda^{-1}]=[-\lambda]$ in $\sqcl$, the resulting relation can be written as $-[\lambda]=[-\lambda]$.
Next, consider the transposition
of the last two vectors of a triple; this transposition maps $t_\lambda$ to
$$
\left(
{1\broose0},
{1\broose\lambda},
{0\broose1}
\right),$$
which can be transformed by $\pmatrix{1&-\lambda^{-1}\cr0&1}$ to
$$
\left(
{1\broose0},
{0\broose\lambda},
{-\lambda^{-1}\broose1}
\right)
=
\left(
{1\broose0},
{0\broose1},
{1\broose-\lambda}
\right).$$
Again, we get $-[\lambda]=[-\lambda]$.
\smallskip
\bf Boundary relations. \rm
A generic quadruple of points in $\P$ can be $G$-transformed to 
$$\left({1\broose0},{0\broose1},{1\broose\lambda},{1\broose\mu}
\right),$$
where genericity is equivalent to $\lambda,\mu\in\dotK$, $\lambda\ne\mu$.
The boundary of the corresponding $3$-simplex is the alternating sum of
four triangles---triples obtained from the quadruple by omitting one element.
We calculate the symbols of those triples.

Omit ${1\broose0}$:
$\left[
\left|\matrix{0&1\cr1&\lambda
}\right|\cdot
\left|\matrix{1&1\cr\lambda&\mu
}\right|\cdot
\left|\matrix{1&0\cr\mu&1
}\right|
\right]=
\left[(-1)\cdot(\mu-\lambda)\cdot1\right]=[\lambda-\mu]$.

Omit ${0\broose1}$:
$
\left[
\left|\matrix{1&1\cr0&\lambda
}\right|\cdot
\left|\matrix{1&1\cr\lambda&\mu
}\right|\cdot
\left|\matrix{1&1\cr\mu&0
}\right|
\right]=
\left[\lambda\cdot(\mu-\lambda)\cdot(-\mu)\right]=[\lambda\mu(\lambda-\mu)]$.

Omit ${1\broose\lambda}$: $[\mu]$.

Omit ${1\broose\mu}$: $[\lambda]$.

The relation is:
$$[\lambda-\mu]-[\lambda\mu(\lambda-\mu)]+[\mu]-[\lambda]=0.$$
Putting $a=\lambda-\mu$ and $b=\mu$ we may rewrite this as
$$[a]+[b]=[a+b]+[ab(a+b)].$$
The relation holds for all $a,b\in\dotK$ that satisfy $a+b\ne0$. This set of relations, plus
the alternation relation $[-a]=-[a]$, gives the classical description of the Witt group $W(K)$ (cf.~[EKM, Theorem 4.8]).
\qed
\smallskip
\bf Definition \a.\t.\rm\par
The tautological second cohomology class of the group $G=PSL(2,K)$ with coefficients in
$U_2=W(K)$ associated (as in Definition 2.1) to the action of $G$ on $X$
will be called the Witt class and denoted by $w$:
$w\in H^2(PSL(2,K),W(K))$.
\smallskip
\bf Remark \a.\t.\rm\par
\item{1)} Let $T$ be the tautological ($W(K)$-valued) cocycle associated to
the $G$-action on $X$. From the proof of Theorem \a.1 it is useful to extract the following
explicit formula for the value of $T$ on a $2$-simplex in $X$ determined by a
triple of pairwise linearly independent vectors $u,v,w\in K^2$:
$$T([u],[v],[w])=[\det(u,v)\det(v,w)\det(w,u)].\leqno(\a.\f)$$
\item{2)} One can see from the proof of Theorem \a.1
that the ordered coefficient
group $U_2^o$ is isomorphic to the Witt--Grothendieck group $\widehat{W}(K)$ of the field $K$.
\item{3)} The space $\P=\P^1(K)$ and the complex $X$ are acted upon by the larger group $PGL(2,K)$.
As a result, the Witt class can be interpreted as a twisted cohomology class of $PGL(2,K)$.
The twisting action of $PGL(2,K)$ on $W(K)$ is easy to see from the formula for $\lambda$ in Lemma \a.2:
the class of $g\in GL(2,K)$ acts on the symbol $[\lambda]$ mapping it to $[\det(g)\cdot\lambda]$. 
\item{4)} For $K=\R$ we have $W(\R)\simeq\Z$. The isomorphism, called the signature map,
maps $[\lambda]$ to $+1$ for $\lambda>0$ and to $-1$ for $\lambda<0$. The pull-back of the Witt class
to $SL(2,\R)$ is a class in $H^2(SL(2,\R),\Z)$;
we will relate it to the usual (topological) Euler class
(see Theorem 13.4 and Fact 13.5).
\item{5)} For $K=\Q$ the Witt group has a large torsion part which is a free summand.
Computer calculations (using the computer algebra
system FriCAS) indicate that the corresponding part of the Witt class is non-trivial.
\smallskip
It is possible to give an explicit formula for a cocycle representing the Witt class.
We use the standard homogeneous resolution (cf.~[Brown, II, \S3]) to describe
group cohomology; $W(K)$-valued $2$-cocycles are then represented by functions
$G\times G\times G\to W(K)$.
\smallskip
\bf Theorem \a.\t.\sl\par
Let us fix an arbitrary non-zero vector $u\in K^2$. The map
$$G\times G\times G\ni(g_0,g_1,g_2)\mapsto
[\det(\tilde{g}_0u,\tilde{g}_1u)\det(\tilde{g}_1u,\tilde{g}_2u)\det(\tilde{g}_2u,\tilde{g}_0u)]\in W(K)\leqno(\a.\f)$$
is a cocycle representing the Witt class $w\in H^2(PSL(2,K),W(K))$. (By $\tilde{g}_i$ we denote an arbitrary lift
of $g_i\in PSL(2,K)$ to $SL(2,K)$. The senseless symbol $[0]$ is interpreted as $0$.)\rm
\smallskip
Proof. It is straightforward to check that the maps
$$\Psi_n\colon(g_0,\cdots,g_n)\mapsto
\cases{ ([\tilde{g}_0u],\ldots,[\tilde{g}_nu]) & if the points $[\tilde{g}_iu]$ are pairwise different,\cr
0& otherwise,}\leqno(\a.\f)
$$
defines a $G$-chain map from the homogeneous standard resolution of $G$ to $C_*X$.
(The only subtle case is when $[\tilde{g}_iu]=[\tilde{g}_ju]$ for exactly one pair of indices $i,j$.
Then $\Psi_{n-1}\partial(g_0,\ldots,g_n)$ has two non-zero summands---however, these summands
cancel in the alternating chain complex $C_*X$.)
Composing $\Psi_2$ with the tautological cocycle $T$ given by (\a.2) we obtain the statement
of the theorem.\qed
\bigskip
\centerline{\bf III. Euler class for ordered fields.}
\bigskip
In this part we define and investigate Euler classes for general linear and
projective groups over arbitrary ordered fields.
\medskip
\bf 8. Tautological Euler classes. Computation of coefficients.\rm
\medskip
\def\a{8}
\m=1
\n=1
Let $K$ be an ordered field.
Let $\G=GL(n,K)$; we will also consider the following closely related
groups (where we put $\dotK=K\setminus\{0\}$, $K_+=\{\lambda\in K\mid \lambda>0\}$):
$$\eqalign{
\PG&:=PGL(n,K)=G/\{\lambda I\mid \lambda\in\dotK\},\cr 
\PpG&:=P_+GL(n,K)=G/\{\lambda I\mid \lambda\in K_+\},\cr 
\Gp&:=GL_+(n,K)=\{g\in \G\mid\det{g}>0\},\cr
\PGp&:=PGL_+(n,K)=\Gp/\{\lambda I\in\Gp\mid \lambda\in K\},\cr
\PpGp&:=P_+GL_+(n,K)=\G_+/\{\lambda I\mid \lambda\in K_+\}.
}
$$
The natural maps between these groups are summarized in the diagram:
$$
\matrix{
\Gp&\longrightarrow&\PpGp&\longrightarrow&\PGp\cr
\downarrow&&\downarrow&&\downarrow\cr
\G&\longrightarrow&\PpG&\longrightarrow&\PG.}
$$
The defining action of $\G$ on $K^n$ induces actions of $\PpG$ on $\Pp$
and of $\PG$ on $\P$; here
$$
\Pp:=\P_+^{n-1}(K)=(K^n\setminus\{0\})/K_+,
\qquad
\P:=\P^{n-1}(K)=(K^n\setminus\{0\})/\dotK,
$$
where the multiplicative groups $K_+$, $\dotK$ act on $K^n$ by homotheties.

Next we define simplicial complexes $X$ and $X_+$ by spanning simplices on
generic tuples of points in $\P$ and $\Pp$, respectively.
We call a tuple $([v_0],\ldots,[v_k])$
generic if every subsequence
of $(v_0,\ldots,v_k)$ of length $\le n$ is linearly independent.
Ordered fields are infinite, therefore these complexes have the star-property
and are contractible.
The complex $X$ is acted upon by $\PG$, and $X_+$ by $\PpG$.
We restrict these actions to $\PGp$ and to $\PpGp$ and we apply the formalism of Part I.
We put
$$U:=U_n(X),\qquad U_+=U_n(X_+)$$
(see Definition 2.1, Definition 1.3),
and we define the Euler classes as tautological classes:
$$\eu:=\tau^n_{ \PGp ,X}\in H^n(\PGp,U),\qquad\eu_+:=\tau^n_{\PpGp,X_+}\in H^n(\PpGp,U_+).\leqno(\a.\f)$$
Notice that $\PGp$ is a normal subgroup of $\PG$. Therefore the group $U$
carries the structure of a $\PG$-module, and the class $\eu$ can also be considered as a
twisted $\PG$-class. (See the second remark after Definition 1.3.) Similarly, the class $\eu_+$ can be
regarded as a twisted $\PpG$-class.

Our first goal is to compute $U$ and $U_+$, as abelian groups and as $\PG$-
and $\PpG$-modules.
\medskip
\bf Theorem \a.\t.\sl\par
Let $U$ ($U_+$) be the coefficient group associated to the
action of $PGL_+(n,K)$ ($P_+GL_+(n,K)$) on the complex of generic tuples
of points in $\P^{n-1}(K)$ ($\P_+^{n-1}(K)$). We have
$$
{U\simeq\cases{0,& for $n$ odd\cr \Z,& for $n$ even}},
\qquad U_+\simeq\Z^{\lfloor{n/2}\rfloor+1}.$$
The $PGL(n,K)$- and $P_+GL(n,K)$-structures are given by
$$[g]\cdot u=\cases{u& if $\det{g}>0$\cr -u& if $\det{g}<0$},\qquad (g\in GL(n,K)).$$
\rm\par
Proof. Both  calculations are based on Fact 2.4. We denote by $(e_1,\ldots,e_n)$ the standard basis of $K^n$.
\smallskip
\bf Calculation of $U$. \rm Non-degenerate ordered simplices of $X$ correspond to generic tuples of points in $\P$.
\smallskip
\bf Lemma \a.\t.\sl\par
The action of $\PGp$ on the set of generic $(n+1)$-tuples of points in $\P$ has one orbit
for $n$ odd and two orbits for $n$ even.
\smallskip\rm
Proof. Let $p=(p_1,\ldots,p_{n+1})$ be a generic $(n+1)$-tuple of points in $\P$. There is an element
$g\in G$ (unique up to scaling) that maps $p$ to the standard tuple $e=([e_1],\ldots,[e_n],[\sum_{i=1}^ne_i])$.
If $n$ is odd we have $\det(-g)=-\det(g)$, therefore $g$ may be chosen in $\Gp$. It follows that
in this case $\PGp$ acts transitively on the set of generic $(n+1)$-tuples.
For $n$ even, all elements $g$ mapping $p$ to $e$ have determinants of the same sign. This sign is a
$\PGp$-invariant of $p$ that we call the sign of $p$ and denote $\sgn(p)$. Generic tuples
of the same sign are $\PGp$-equivalent: if $\sgn(p)=\sgn(p')$, $gp=e$, $gp'=e$, then
$g^{-1}g'p'=p$ and $\sgn\det(g^{-1}g')=+1$.
\phantom{aaaaa}\hfill\qed{(Lemma)}
\smallskip
The case of $n$ odd is now straightforward. The image of any $n$-simplex of $X$ in $(C_nX)_{\PGp}$ is one and the same
generator of that cyclic group. The boundary of an $(n+1)$-simplex of $X$ is an alternating sum of
an odd number of $n$-simplices, hence its image in $(C_nX)_{\PGp}$ is again that generator. It follows
that $U=0$ for $n$ odd.

Suppose now that $n$ is even. Lemma 8.2 and Fact 2.4 imply that $U$ is the quotient of the free abelian group
with two generators by two sets of relations (boundary relations and alternation relations).
The generators correspond to (representatives of) $\PGp$-orbits on the set of generic $(n+1)$-tuples
of points in $\P$; explicitly, we choose
$$([e_1],\ldots, [e_n],[\sum_{i=1}^ne_i])\leqno(\a.\f)$$
and denote it $[+]$ or $[+1]$, and
$$([e_1],\ldots,[e_{n-1}], [-e_n],[(\sum_{i=1}^{n-1}e_i)-e_n])\leqno(\a.\f)$$
and denote it $[-]$ or $[-1]$. We call $[+]$ and $[-]$ symbols.

The group $\Z^2$ generated by $[+]$ and $[-]$ is isomorphic to $(C^o_nX)_{\PGp}$. The image in this group
of an ordered $n$-simplex of $X$ corresponding to a generic $(n+1)$-tuple $p=(p_1,\ldots,p_{n+1})$
of points in $\P$ is $[\sgn(p)]$. In practice, the sign can be calculated as follows:
let $p_i=[v_i]$ for $v_i\in K^n$, and let $v_{n+1}=\sum_{i=1}^n\alpha_iv_i$; then
$$\sgn(p)=\sgn(\det(v_1,\ldots,v_n)\cdot\alpha_1\cdot\ldots\cdot\alpha_n)\leqno(\a.\f)$$
(this is the sign 
of the determinant of the matrix mapping $(v_1,\ldots,v_{n+1})$ to $(e_1,\ldots,e_n,\sum_{i=1}^ne_i)$;
that matrix is the inverse of the product of the matrix with columns $(v_1,\ldots,v_n)$
and the diagonal matrix with diagonal entries $(\alpha_1,\ldots,\alpha_n)$).

The alternation relations all reduce to $-[+]=[-]$. Indeed, it is straightforward to check that
transposing two neighbouring elements in $(\a.2)$ or $(\a.3)$ changes the sign of the tuple.

We now discuss the boundary relations. We observe that any ordered non-degenerate $(n+1)$-simplex of
$X$---corresponding to a generic $(n+2)$-tuple of points in $\P$---can be mapped by an element
of $\PGp$ to
$$\Delta=([e_1],[e_2],\dots,[e_{n-1}],[se_n],[(\sum_{i=1}^{n-1}e_i)+se_n],
[(\sum_{i=1}^{n-1}b_ie_i)+sb_ne_n]).$$
Here $s=\pm1$ and $b_i\in K$;
genericity means that $b_i\neq 0$ and  $b_i\neq b_j$.
We have
$$\partial\Delta=\sum_{j=1}^{n+2}(-1)^{j-1}[s_j],\leqno(\a.\f)$$
where $s_j$ is the sign of the tuple obtained from $\Delta$ by omitting
the $j$-th element. We have $s_{n+2}=s$ and $s_{n+1}=s\sgn(\prod b_i)$.
We claim the for $j<n+1$ we have
$$s_j=(-1)^j s \sgn(b_j\prod_{i\ne j} (b_i-b_j)).\leqno(\a.\f)$$
Indeed, for $j<n$ we have
$$\eqalign{
&\sgn\det(e_1,\ldots,\widehat{e_j},\ldots,e_{n-1},se_n,\sum_{i=1}^{n-1} e_i+se_n)
=(-1)^{n-j}s=(-1)^js,\cr
&\sum_{i=1}^{n-1}b_ie_i+sb_ne_n=
b_j(\sum_{i=1}^{n-1} e_i+se_n)+\sum_{i\ne j,n}(b_i-b_j)e_i+(b_n-b_j)se_n,}$$
while for $j=n$ we have 
$$\eqalign{
&\sgn\det(e_1,\ldots,e_{n-1},\widehat{se_n},\sum_{i=1}^{n-1} e_i+se_n)=s=
(-1)^ns,\cr
&\sum_{i=1}^{n-1}b_ie_i+sb_ne_n=
b_n(\sum_{i=1}^{n-1} e_i+se_n)+\sum_{i=1}^{n-1}(b_i-b_n)e_i.}$$

Putting (\a.5) and (\a.6) together we get:
$$
\eqalign{
\partial\Delta=\sum_{j=1}^n(-1)^{j-1}[(-1)^js\sgn{b_j}\prod_{i\ne j}\sgn(b_i-b_j)]+(-1)^n[s\sgn(\prod_ib_i)]-[s].
}\leqno(\a.\f)$$
We will show that this relation is trivial, i.e.~that all the symbols cancel.
We can assume $s=+1$; the case $s=-1$ will automatically follow. Indeed, changing $s$ 
flips all the symbols: $[+1]\leftrightarrow[-1]$, and trivially transforms a trivial relation to
a trivial relation.

Let us artificially put $b_{n+1}=0$; then we may rewrite $(\a.7)$ more uniformly:
$$\partial\Delta=\sum_{j=1}^{n+1}(-1)^{j-1}[\sgn(\prod_{i=1}^{j-1}(b_j-b_i)\prod_{i=j+1}^{n+1}(b_i-b_j))]-[+1].\leqno(\a.\f)$$
Let $\sigma\in S_{n+1}$ be the permutation ordering the numbers (indices)  in the
same way that the sequence $b$ does: $\sigma(i)<\sigma(k)\iff b_i<b_k$.
We put $\inv(j)=\#\{i\mid (i-j)(\sigma(i)-\sigma(j))<0\}$ (the number of inversions
of $\sigma$ in which $j$ is involved). Then our relation is
$$\sum_{j=1}^{n+1}(-1)^{j-1}[(-1)^{\inv(j)}]-[+1].\leqno(\a.\f)$$
\bf Lemma \a.\t. \sl
$(-1)^{\inv(j)}=(-1)^{\sigma(j)-j}$.
\rm\par
Proof. If exactly $k$ of the indices smaller than $j$ are mapped by $\sigma$ to indices
larger than $\sigma(j)$, then $\sigma(j)-j+k$ of the ones larger than $j$
have to be mapped to values smaller than $\sigma(j)$.
Then $\inv(j)=k+\sigma(j)-j+k\equiv\sigma(j)-j\,(\hbox{\rm mod}\,2)$.
\qed{(Lemma)}
\smallskip
Thus, $\inv(j)$ is odd if and only if $j$ and $\sigma(j)$ differ in parity.
Since the number of even $j$'s is equal to  the number of even $\sigma(j)$'s,
this difference in parity appears equally often in each of the two  forms:
$(j,\sigma(j))=({\rm odd},{\rm even})$, $(j,\sigma(j))=({\rm even},{\rm odd})$.
In (\a.9), pairs of the first kind lead to summands $+[-]$, while pairs of the second kind
give $-[-]$. Thus, all the appearances of the symbol $[-]$ cancel. It follows
that the sum adds up to $[+]$, which is cancelled by the extra term.

We have shown that the boundary relations are trivial. It follows that $U$ is the quotient
of $\Z^2$ by the alternation relations, i.e.~$U\simeq\Z$.
The $\PG$-structure description follows from the formula
$\sgn(gp)=\sgn(\det{g})\cdot\sgn(p)$, valid for $g\in G$ and all generic
$(n+1)$-tuples $p$ of points in $\P$.
\medskip
\bf Calculation of $U_+$. \rm 
A generic $n$-tuple of points in $\Pp$ can be lifted to 
$n$ linearly independent vectors in $K^n$. The matrix $M$ with columns given by these vectors is well-defined
up to multiplication on the right by diagonal matrices with positive entries on the diagonal.
The sign of $\det{M}$ is thus an invariant of the tuple (the \it sign of the tuple\rm); it is also a $\Gp$-invariant of the tuple
($\det{gM}=\det{g}\cdot\det{M}=\det{M}$ for $g\in\Gp$). In fact, this sign is the full $\Gp$-invariant:
the tuple of vectors is transformed to the standard basis by $M^{-1}$ if $\det{M}>0$,
and to the basis $(e_1,\ldots,e_{n-1},-e_n)$ by $M^{-1}$ with negated lowest row if $\det{M}<0$.
We have shown the following statement.
\smallskip
\bf Lemma \a.\t.\sl\par
The action of $\Gp$ on the set of generic $n$-tuples of points in $\Pp$ has two orbits, 
detected by the sign of the tuple.\rm
\smallskip
We now consider the $\Gp$-action on the set of generic $(n+1)$-tuples.
The symbol of such a tuple $(p_1,\ldots,p_{n+1})$ is defined to be a sequence of $n+1$ signs:
$[s;s_1,\ldots,s_n]$. Here $s$ is the sign of $(p_1,\ldots,p_n)$. To get the other signs,
we first lift each $p_i$ to a vector $v_i$. Then we express $v_{n+1}$ in terms of the other $v_i$:
$v_{n+1}=\sum_{i=1}^n a_nv_n$. Finally, $s_i=\sgn{a_i}$ (genericity implies $a_i\ne0$). 
Clearly, the symbol is $\Gp$-invariant.

\smallskip
\bf Lemma \a.\t. \sl\par
Two generic $(n+1)$-tuples of points in $\Pp$ are $\Gp$-equivalent
if and only if they have the same symbol.
\rm\par
Proof. 
Let the tuples be $(p_i)$ and $(q_i)$, with symbol $[s;s_1,\ldots,s_n]$.
Then a lift of $(p_i)$ is equivalent to $(e_1,\ldots,e_{n-1},se_n,\sum_{i=1}^na_ie_i)$
for some $a_i\in \dotK$, while a lift of $(q_i)$ is equivalent
to $(e_1,\ldots,e_{n-1},se_n,\sum_{i=1}^nb_ie_i)$
for some $b_i\in \dotK$; moreover, $s_i$ is the (common) sign of $a_i$ and of $b_i$
(for $i<n$; and $s_n$ is the common sign of $sa_n$ and $sb_n$).
These two representing tuples of vectors are projectively related
by the diagonal matrix with positive diagonal entries $b_i/a_i$.
\qed{(Lemma)}
\smallskip
The following observation describes the $\G$-action on symbols and allows us
to determine the $\PpG$-structure on $U_+$.
\smallskip
\bf Fact \a.\t. \sl\par
The symbol of a generic $(n+1)$-tuple is $\G$-equivariant:
if $g\in\G$, $\det{g}<0$, then the tuples $(p_i)$ and
$(gp_i)$ have the opposite leading symbol sign $s$, and coinciding
remaining symbol signs.
\rm
\smallskip
It follows from Lemma \a.5 that the group $U_+$ is the quotient of the free abelian group
spanned by symbols by alternation and boundary relations. We first deal with the alternation relations.

\bf Alternation relations. \rm The symbol $[s;s_1,\ldots,s_n]$ is represented by the tuple
$$(e_1,\ldots,se_n,\sum_{i=1}^{n-1}s_ie_i+ss_ne_n).$$
Suppose that this tuple is permuted; what happens to the symbol?
Since permutation commutes with `linear map applied to each element', we may and will assume $s=+1$---in our arguments,
but not on the final statements.
We first treat the case of a permutation $\sigma$ that fixes the last element.
Then the symbol of the permuted tuple is $[\sgn{\sigma};(s_{\sigma^{-1}(i)})]$.
We get the `usual permutation relation' 
$$[s;(s_i)]=\sgn{\sigma}[s\sgn{\sigma};(s_{\sigma^{-1}(i)})].$$
Now let us consider the transposition of $k$ and $n+1$.
The new leading sign is $$\det(e_1,\ldots,e_{k-1},\sum_{i=1}^ns_ie_i,e_{k+1},\ldots,e_n)=s_k.$$ 
We also have
$$e_k= s_k\sum_{i=1}^ns_ie_i+
\sum_{i\ne k}(-s_ks_i)e_i,$$  
so that the total symbol after transposition is
$$[s_k; -s_ks_1,\ldots,-s_ks_{k-1},s_k,-s_ks_{k+1},\ldots,-s_ks_n].$$
The `transposition relation' is thus
$$[s;(s_i)]=-[ss_k; -s_ks_1,\ldots,-s_ks_{k-1},s_k,-s_ks_{k+1},\ldots,-s_ks_n].$$
In words: if the $k$'th sign is $+$, then we can flip all the other signs 
(except the leading sign); the resulting symbol will be equal to minus the original.
If the $k$'th sign is $-$, we get $[s;(s_i)]=-[-s;(s_i)]$.

There is a difference between the cases $n=2$ and $n>2$. In the latter case,
for any sequence of $n$ signs there exists a stabilizing transposition; therefore,
any sequence of $n$ signs can be ordered (put in the form $++\ldots--$) by
an even permutation. Let us begin with the case $n>2$.
\smallskip
\bf The case $n>2$. \rm
As already mentioned, in this case one can use the usual permutation relation
to order the non-leading signs of a symbol without changing the leading sign.
To shorten the notation, we will use $a^{+}$ for $[+;++\ldots--]$
($a$ plus signs after the semi-colon),
and $a^{-}$ for $[-;++\ldots--]$ ($a$ plus signs after the semi-colon).
For example, when $n=3$, we put $0^+=[+;---]$, $2^+=[+;++-]$ and $2^-=[-;++-]$.
The transposition relation (with $s_k=+1$) gives $a^\pm=-(n-a+1)^\pm$ (for $a>0$).
Picking $s_k=-1$ in the transposition relation we get $a^+=-a^-$---for $a<n$, but
$n^+=-n^-$ also holds, due to $n^+=-1^+=1^-=-n^-$.  
To summarize:
\smallskip
\bf Lemma \a.\t. \sl\par
Let $n>2$. Let $A=\{a^+\mid 0\le a\le \lfloor{n/2}\rfloor\}$.
The quotient of the group $(C^o_nX_+)_{\PpGp}$ by the set of alternation relations
is the free abelian group with generating set $A$ for $n$ even; 
for $n$ odd it is the direct sum of the free abelian group generated by $A$ and
a $\Z/2$ generated by $\left({n+1\over2}\right)^+$.\rm
\smallskip
(The extra $\Z/2$-summand appearing for $n$ odd will eventually get killed
by the boundary relations). 
\smallskip
\bf The case $n=2$. \rm
There are eight symbols. The usual permutation relation for $\sigma=(12)$ gives:
$$
\eqalign{
[+;++]=-[-;++],
\qquad
[-;-+]&=-[+;+-],
\qquad
[-;+-]=-[+;-+],
\cr[+;--]&=-[-;--].}
$$
The transposition relation (for $\sigma=(23)$) is
$[s;s_1,s_2]=-[ss_2;-s_2s_1,s_2]$. This gives:
$$
\eqalign{
-[-;++]=[-;-+],
\qquad
-[+;+-]&=[-;+-],
\qquad
-[+;-+]=[+;++],
\cr[+;--]&=-[-;--].
}
$$
We see that all the (signed) symbols appearing in the first lines
of the above formulae are identified. In particular, $[+,+-]=[+;-+]$, so 
that the $a^\pm$ notation still makes sense. Also, the relations
$a^++a^-=0$, $1^++2^+=0$ can be read off from the ones displayed above.
Therefore the conclusion of Lemma \a.7 holds for $n=2$.
\smallskip
\bf Boundary relations. \rm We will show that they all follow from the alternation relations
(with the exception of $\left({n+1\over2}\right)^+=0$).
Let us calculate the boundary of an $(n+1)$-simplex of $X_+$ represented by a generic $(n+2)$-tuple of vectors.
Such a tuple of vectors can be transformed by an element of $\Gp$
to $\Delta=(e_1,\ldots,e_{n-1},se_n,\sum_{i=1}^ns_ie_i,\sum_{i=1}^ns_ib_ie_i)$.
The genericity condition (assuming $s=+1$) is: all $b_i$ non-zero and pairwise different.
(This will follow from the calculation of $\partial\Delta$.)
If $s=-1$, we can transform the tuple by an orientation changing linear map;
this will change all leading signs in $\partial\Delta$, and not touch the other signs.
Thus, we will assume $s=+1$---and then double the set of the resulting
boundaries by changing the leading signs.
If we omit $e_j$ ($i\le n$) from $\Delta$, then the sign of the determinant
of the standardizing matrix is the same as that of 
$$\det(e_1,\ldots,e_{j-1},e_{j+1},\ldots,e_n, \sum_{i=1}^ns_ie_i)=(-1)^{n-j}s_j.$$
The other signs can be read off from:
$$\sum_{i=1}^n s_ib_ie_i=\sum_{i\ne j}(s_ib_i-s_ib_j)e_i+b_j\sum_{i=1}^ns_ie_i.$$
The total symbol (for $j$ omitted) is thus:
$[(-1)^{n-j}s_j;(s_i\sgn(b_i-b_j))_{i\ne j},\sgn{b_j}]$.

Omitting the $(n+1)$'st element gives $[+1;(s_i\sgn{b_i})]$.

Omitting the $(n+2)$'nd element yields $[+1;(s_i)]$.
So, finally
$$
\eqalign{
\partial\Delta=\sum_{j=1}^n(-1)^{j-1}&[(-1)^{n-j}s_j;(s_i\sgn(b_i-b_j))_{i\ne j},\sgn{b_j}]\cr
&+(-1)^n[+1;(s_i\sgn{b_i})]+(-1)^{n+1}[+1;(s_i)].
}$$
Let us rewrite this boundary relation. We put (artificially) $b_{n+1}=0$ and
$s_{n+1}=-1$. Then we have $\sgn{b_j}=s_{n+1}\sgn(b_{n+1}-b_j)$, $s_i\sgn{b_i}=s_i\sgn(b_i-b_{n+1})$ and
$(-1)^{n-(n+1)}s_{n+1}=+1$,
so that 
$$
\eqalign{
\partial\Delta
&=\sum_{j=1}^{n+1}(-1)^{j-1}[(-1)^{n-j}s_j;(s_i\sgn(b_i-b_j))_{i\ne j}]+(-1)^{n+1}[+1;(s_i)]\cr
&=\sum_{j=1}^{n+1}(-1)^{n-1}[s_j;(s_i\sgn(b_i-b_j))_{i\ne j}]+(-1)^{n+1}[+1;(s_i)],
}
$$
where the last equality uses $a^+=-a^-$. Thus, we need to deduce 
from our alternation relations (assuming $s_{n+1}=-1$, $b_{n+1}=0$) that 
$$\sum_{j=1}^{n+1}[s_j;(s_i\sgn(b_i-b_j))_{i\ne j}]+[+1;(s_i)_{i\ne n+1}]=0.$$
We see that we may change the order of summation to follow the
order of the $b_i$. Indeed, the above summation can be phrased in an index-free way as
follows. We have a set of $n+1$ numbers, each with an attached sign
(one of these pairs being $0$ with $-$). For each element $x$ of the set
we form the corresponding symbol $k^s$, where $s$ is the sign attached to the element $x$,
and $k$ is the number of positive expressions of the form $t(y-x)$, where $y$ runs through
our set (and $y\ne x$) and $t$ is the sign attached to $y$. Finally, there is an extra summand
$\ell^+$ with $\ell$ counting all the positive signs.

We may thus renumber the $b_i$ and the $s_i$ (in the same way) so as to have
$b_1>b_2>\ldots>b_n>b_{n+1}$, with an unknown $b$ equal to zero
and the corresponding $s$ equal to $-1$ and omitted in the extra summand $[+1;(s_i)]$.
Let us now consider two consecutive summands (numbered $j$ and $j+1$).
They differ at most by the leading sign: $s_j$ versus $s_{j+1}$, and by the $j$'th non-leading sign:
$s_{j+1}\sgn(b_{j+1}-b_j)=-s_{j+1}$ versus $s_j\sgn(b_j-b_{j+1})=s_j$.
Substituting all four possible combinations of $(s_j,s_{j+1})$ we get:
\smallskip
\bf Claim. \sl Two consecutive summands are of one of the forms:
$$(a^\pm,a^\mp),\quad (a^+,(a+1)^+),\quad(a^-,(a-1)^-).$$
\smallskip\rm
Suppose that $k$ of the $s_i$'s are positive. Then
the extra summand is $k^+$, while the first and the last one depend on $(s_1,s_{n+1})$ and are:
$$\matrix{
s_1&s_{n+1}&{\rm first}&{\rm last}\cr
+&+&(n-k+1)^+&(k-1)^+\cr
+&-&(n-k+1)^+&k^-\cr
-&+&(n-k)^-&(k-1)^+\cr
-&-&(n-k)^-&k^-
}
$$
We can append the extra summand $k^+$ to the sum (while keeping the rule of the claim)
and get summation starting from $(n-k+1)^+$ or $(n-k)^-$ and ending at $k^+$.
Then we start cancelling consecutive pairs $(a^\pm,a^\mp)$ (except that we do not
cancel the first and the last element) until the sequence becomes monotone
(possibly except the first or the last pair).
If the final monotone sequence runs from
$(n-k+1)^+$ to $k^+$ then the terms pairwise cancel (first with last, second with last-but-one, etc.)
if $n$ is even, and ${n+1\over2}$ is left if $n$ is odd.
If the sequence starts with $(n-k)^-$, we may put an extra pair $((n-k+1)^+,(n-k+1)^-)$ at the beginning
of the sequence, to reduce to the former case---except when $k=0$. 
If $k=0$, we get a sequence running from $n^-$ to $0^+$, ie. $(n^-,(n-1)^-,\ldots,1^-,0^-,0^+)$.
The first $n$ terms cancel in the same manner as before, and $0^-+0^+=0$.

Finally, since the set of permutation relations is invariant under the `exponent sign' flip ($a^{\pm}\leftrightarrow a^{\mp}$),
the boundary relations obtained from $\Delta$ with $s=-1$ are dealt with in the same way.

Fact \a.6 and the relation $a^-=-a^+$ imply the remaining claim of the theorem (the one describing
the $\PpG$-structure on $U_+$).
\qed({Theorem \a.1)}
\smallskip
\bf Remark \a.\t.\rm\par
Let $T$ be the tautological ($U$-valued) $n$-cocycle associated to the $\PGp$-action on $X$,
and let $T_+$ be the tautological ($U_+$-valued) $n$-cocycle associated to the $\PpGp$-action on $X_+$.
From the proof of Theorem \a.1 it is useful to extract the following explicit description of these cocycles.
\item{a)} Let $n$ be even; then $U\simeq \Z$ is generated by the symbol $[+]$.
Suppose that $\sigma=([v_1],\ldots,[v_{n+1}])$ is an $n$-simplex of $X$. Then
$v_{n+1}=\sum_{i=1}^n\alpha_iv_i$ for some $\alpha_i\in\dotK$. We have (see (\a.4)):
$$T(\sigma)=[\sgn(\det(v_1,\ldots,v_n)\cdot\alpha_1\cdot\ldots\cdot\alpha_n)].$$
(Recall that $[-]=-[+]$.)
\item{b)} Recall that $U_+\simeq\Z^{\lfloor{n/2}\rfloor+1}$ with free generating set
$A=\{a^+\mid a=0,\ldots,\lfloor{n/2}\rfloor\}$ (see Lemma \a.7).
Suppose that $\sigma=([v_1],\ldots,[v_{n+1}])$ is an $n$-simplex of $X_+$. Then
$v_{n+1}=\sum_{i=1}^n\alpha_iv_i$ for some $\alpha_i\in\dotK$. To $\sigma$ we assign an
$(n+1)$-tuple of signs $[s;s_1,\ldots,s_n]$, where
$s=\sgn{\det(v_1,\ldots,v_n)}$ and 
$s_i=\sgn{\alpha_i}$. Next we put
$T_+(\sigma)=a^+$ (if $s=+1$ and $a$ of the $s_i$ are $+1$)
or $T_+(\sigma)=a^-$ (if $s=-1$ and $a$ of the $s_i$ are $+1$). Finally, we express
the symbol in term of the elements of $A$ using the relations:
$a^+=-a^-$; $a^\pm=-(n+1-a)^\pm$ (for $a>0$).
\smallskip
\bf Definition \a.\t.\rm\par
The splitting of $U_+=\oplus_{k=0}^{\lfloor{n/2}\rfloor}\Z k^+$
into cyclic summands generated by the elements $k^+$ ($0\le k\le\lfloor{n/2}\rfloor$)
induces the corresponding splittings of cocycles and cohomology classes:
$$\eqalign{T_+&=\oplus T_k,\quad T_k\in Z^n(\Hom{}_{\PpGp}(C_*X_+,\Z));\cr
\eu_+&=\oplus \eu_k,\quad \eu_k\in H^n(\PpGp,\Z);\cr
\tau_+&=\oplus \tau_k,\quad \tau_k\in H^n(\Hom{}_{\PpGp}(C_*X_+,\Z)).}$$
(In the last formula $\tau_+$ ($\tau_k$) is the cohomology class
of $T_+$ ($T_k$).)
\smallskip
\bf Remark \a.\t.\rm\par
Suppose that $K<L$ is a field extension, and that on $K$ and on $L$ there are compatible field orders.
Then we have the group embedding $\phi\colon\PpGp(K)\to\PpGp(L)$, and
the natural $\phi$-equivariant simplicial complex embedding $X_+(K)\to X_+(L)$ inducing
a coefficient group map $f\colon U_+(K)\to U_+(L)$. But, in our field--independent description of $U_+$
(see Theorem \a.1 and Remark \a.8) the map $f$ is represented by identity.
Applying Theorem 1.5 to these data
we obtain $\phi^*\eu_+(L)=\eu_+(K)$---the Euler class $\eu_+$ is stable under ordered field restriction.
It follows that all $\eu_k$ are also stable. Analogous arguments show the same stability statement
for $\eu$.
\smallskip
\bf Remark \a.\t.\rm\par
It follows from Theorem 3.1 that the classes $\eu$, $\eu_+$ and $\eu_k$ are
bounded.
\medskip
\bf 9. Relation between the classes $\eu_k$.\rm
\medskip
\def\a{9}
\m=1
\n=1
The classes $\eu_k$ defined in Definition 8.9 are related.
\smallskip
\bf Theorem \a.\t.\sl\par
$$\sum_{k=0}^{\lfloor{n/2}\rfloor}(n-2k+1)\eu_k=0$$
\rm\par
Proof.
We will see that this relation holds already in  $H^n(\Hom{}_{\PpGp}(C_*X_+,\Z))$
for the classes $\tau_k$.
To prove it, we  will find a cochain $c\in C^{n-1}(\Hom{}_{\PpGp}(C_*X_+,\Z))$
such that
$$\delta c=\sum_{k=0}^{\lfloor{n/2}\rfloor}(n-2k+1)T_k\leqno(\a.\f)$$
in $C^n(\Hom{}_{\PpGp}(C_*X_+,\Z))$.
The boundary map
$$\partial\colon (C_nX_+)_{\PpGp}\to (C_{n-1}X_+)_{\PpGp}$$
factors as the composition of the projection $(C_nX_+)_{\PpGp}\to(C_nX_+/B_nX_+)_{\PpGp}=U_+$
and a map 
$\partial'\colon U_+\to (C_{n-1}X_+)_{\PpGp}$. Each $T_k$ also factors---as the composition of the same projection
and the projection $T_k'$ of $U_+$ on the $k^+$-summand.
Recall that $(C_{n-1}X_+)_\PpGp\cong\Z$ (with generator $[+]$; see Lemma 8.4).
We now consider a generator $a^+$ of $U_+$ and determine $\partial'(a^+)$.
Let $v_a=e_1+\ldots+e_a-(e_{a+1}+\ldots+e_n)$; then $(e_1,\ldots,e_n,v_a)$ determines a simplex in $X_+$
representing $a^+$. We have
$$\eqalign{
\partial'(a^+)
&=\partial[e_1,\ldots,e_n,v_a]\cr
&=\sum_{j=1}^n(-1)^{j+1}[e_1,\ldots,\widehat{e_j},\ldots,e_n,v_a]+(-1)^n[e_1,\ldots,e_n]\cr
&=\sum_{j=1}^n(-1)^{j+1}(-1)^{n-j}[e_1,\ldots,v_a,\ldots,e_n]+(-1)^n[+]\cr
&=\sum_{j=1}^a(-1)^{n+1}[+]+\sum_{j=a+1}^n(-1)^{n+1}[-]+(-1)^n[+]\cr
&=(-1)^n((1-a)[+]-(n-a)[-])=(-1)^n(n-2a+1)[+]
}
$$
Let $c\in C^{n-1}(\Hom{}_\PpGp(C_*X_+,\Z))={\rm Hom}((C_{n-1}X_+)_\PpGp,\Z)$ be defined by
$c([+])=(-1)^n$. Then $(c\circ\partial')(a^+)=c(\partial' a^+)=(n-2a+1)=\sum_{k=0}^{\lfloor{n/2}\rfloor}(n-2k+1)T_k'(a^+)$
holds for each $a^+$. Formula (\a.1) follows.\qed
\medskip
\bf 10. The Smillie argument.\rm
\medskip
\def\a{10}
\m=1
\n=1
The Smillie argument (see [Gro, Section 1.3]) can be used to show that the classes $\eu_k$ are proportional
in a weak sense.
\smallskip
\bf Theorem \a.\t.\sl\par
For any $h\in H_n(B\PpGp,\Z)$ (or $h\in H_n(B\PpGp,\Zm(m))$ for $m$ odd) and any $k\le\lfloor{n/2}\rfloor$
$$\langle \eu_k,h\rangle=(-1)^k{n+1\choose k}\langle \eu_0,h\rangle.$$
If $n$ is odd, then $\langle \eu_k,h\rangle=0$ for all $k$.
\smallskip
\rm Proof. It is well-known that there exists a finite simplicial complex
$Y$, a simplicial cycle $Z\in Z_n(Y,\Z)$ (or in $Z_n(Y,\Zm(m))$), and a map $f\colon Y\to BG$, such that
$f_*[Z]=h$. Let $P=f^*E\PpGp$ (the pull-back of the universal bundle over $B\PpGp$).
Then
$$\langle \eu_k,h\rangle=\langle \eu_k,f_*[Z]\rangle=
\langle f^*\eu_k,[Z]\rangle=\langle \eu_k(P),[Z]\rangle$$
for each $k$. We will use Theorem 4.4 to compute
$\langle \eu_k(P),[Z]\rangle$. Let $E=P\times_\PpGp\P_+$ be the associated
bundle.

Pick a generic section $s\colon Y^{(0)}\to E$. Genericity means that the values of the section
at the vertices of any $n$-simplex of $Y$, viewed as points in $\P_+$ via a flat trivialization
of $E$ over that simplex, form a generic tuple of points. Such a section can be picked vertex-by-vertex.
At a vertex $y$ the genericity conditions mean that a certain finite union of proper projective subspaces
of $E_y$ is prohibited; since the ordered field $K$ is infinite, that union does not fill out $E_y$
and a generic choice is possible. Any generic section $s$ determines a simplicial section
of the associated $X_+$-bundle over $Y$, and then Theorem 4.4 may be applied.

For any
function $\epsilon\colon Y^{(0)}\to\{\pm1\}$ we can form a modified section
$\epsilon s\colon Y^{(0)}\to E$ defined in the obvious way: if $s(v)=(p,[q])$ (for some $q\in K^n\setminus\{0\}$),
then $(\epsilon s)(v)=(p,[\epsilon(v)q])$. Every section $\epsilon s$ is again generic.
Theorem 4.4 gives $\eu_+(P)=[(\epsilon s)^*T_+]$, and coefficient splitting allows us to deduce
$\eu_k(P)=[(\epsilon s)^*T_k]$; both formulae hold for all functions $\epsilon$.
For a given $n$-simplex $\sigma$ of $X_+$ we will
average the expression $\langle(\epsilon s)^*T_k,\sigma\rangle$ over all
possible functions $\epsilon$.

Let $\sigma=(v_1,\ldots,v_n,v_{n+1})$ be one of the $n$-simplices of $Y$.
Let us choose a flat section $r$ of $P$ over $\sigma$, and let $s(v_i)=[r,[q_i]]$,
for $q_i\in K^n\setminus\{0\}$.
We denote by $s_*\sigma$ the $n$-simplex of $X_+$ given by
$([q_1],\ldots,[q_n],[q_{n+1}])$. This definition depends on the choice of $r$,
but different choices lead to simplices equivalent under the
$\PpGp$-action.
The expression $\langle T_k,s_*\sigma\rangle$ is well-defined and equal to
$\langle s^*T_k,\sigma\rangle$.
Let $\eta=\sgn\det(q_1,\ldots,q_n)$, and let
$q_{n+1}=\sum_{i\le n}a_iq_i$. Suppose that exactly $\ell$ of the coefficients
$a_i$ are positive---so that the symbol of $s_*\sigma$ is $\ell^\eta$.

For any function $\epsilon$ we have
 $(\epsilon s)_*\sigma=([\epsilon_1q_1],\ldots,[\epsilon_nq_n],[\epsilon_{n+1}q_{n+1}])$,
where $\epsilon_i=\epsilon(v_i)$.
We wish to determine all functions $\epsilon$ such that
$\langle T_k,(\epsilon s)_*\sigma\rangle\ne0$. This will happen
if and only if the decomposition
$\epsilon_{n+1}q_{n+1}=\sum_{i\le n}b_i\epsilon_iq_i$ has either $k$ or $n+1-k$ positive
coefficients $b_i$.

Let us first focus on the case of $\epsilon_{n+1}=+1$ and $k$ positive $b_i$'s.
We will represent the appropriate functions $\epsilon$ in the form $\epsilon'\epsilon''$;
the idea is that $\epsilon'$ makes all the non-leading signs negative, while $\epsilon''$
changes $k$ of them to $+$. In more detail:
$\epsilon'_{n+1}=+1$ and $\epsilon'_i=-\sgn{a_i}$ for $i\le n$, while
$\epsilon''$ is arbitrary with $k$ negative and $n-k$ positive values
(plus $\epsilon''_{n+1}=+1$). For such $\epsilon'$ and $\epsilon''$ the symbol
of
$(\epsilon'\epsilon''s)_*\sigma$ is $k^\pm$, where the exponent is
$\prod_{i\le n}\epsilon'_i\prod_{i\le n}\epsilon''_i\cdot\sgn\det(q_1,\ldots,q_n)=(-1)^\ell(-1)^k\eta$.
There are ${n\choose k}$ appropriate functions $\epsilon$.

For $\epsilon_{n+1}=+1$ and $n+1-k$ positive $b_i$'s we get ${n\choose n+1-k}={n\choose k-1}$
possibilities, yielding $(n+1-k)^\pm$, with the exponent equal to
$(-1)^\ell(-1)^{n+1-k}\eta=-(-1)^n(-1)^\ell(-1)^k\eta$.

If $\epsilon_{n+1}=-1$ the analysis is analogous. The difference is that
$\epsilon'$ should now be: $\epsilon_i'=\sgn{a_i}$; therefore, the only
change is $(-1)^{n-\ell}$ instead of $(-1)^\ell$ in the final exponent sign formulae.

Putting these together, we get (with $N=\# Y^{(0)}$):
$$\eqalign{
\langle T_k,\sum_\epsilon(\epsilon s)_*\sigma\rangle
=&
\langle T_k,
2^{N-(n+1)}\Biggl(
{n\choose k}k^{(-1)^\ell(-1)^k\eta}+
{n\choose k-1}(n+1-k)^{-(-1)^n(-1)^\ell(-1)^k\eta}\cr
&+{n\choose k}k^{(-1)^{n-\ell}(-1)^k\eta}+
{n\choose k-1}(n+1-k)^{-(-1)^n(-1)^{n-\ell}(-1)^k\eta}
\Biggr)\rangle\cr
=&
2^{N-(n+1)}
\langle T_k,
\Biggl(
{n\choose k}((-1)^\ell(-1)^k\eta+(-1)^{n-\ell}(-1)^k\eta)k^+\cr
&-{n\choose k-1}((-1)^n(-1)^\ell(-1)^k\eta+(-1)^n(-1)^{n-\ell}(-1)^k\eta)(n+1-k)^+
\Biggr)\rangle\cr
=&2^{N-(n+1)}
\langle T_k,
\Biggl(
{n\choose k}(1+(-1)^n)(-1)^\ell(-1)^k\eta k^+\cr
&-{n\choose k-1}(1+(-1)^n)(-1)^n(-1)^\ell(-1)^k\eta(n+1-k)^+
\Biggr)\rangle.
}
$$
For $n$ odd, this is zero. Thus, we assume $n$ even; then the coefficients in the above expression add up to
$$2^{N-n}\left({n\choose k}+{n\choose k-1}\right)(-1)^k(-1)^\ell\eta=
2^{N-n}{n+1\choose k}(-1)^k(-1)^\ell\eta.$$
Similarly, if $n+1-\ell$ of the coefficients $a_i$ are positive, we get
$$2^{N-n}{n+1\choose k}(-1)^k(-1)^{n+1-\ell}\eta.$$
In both cases, the result can be expressed as
$$2^{N-n}{n+1\choose k}(-1)^k\langle (-1)^\ell T_\ell,s_*\sigma\rangle.$$
Since on any $s_*\sigma$ exactly one of $T_\ell$ is non-zero, we
can rewrite this formula as follows:
$$2^{N-n}{n+1\choose k}(-1)^k\sum_\ell\langle (-1)^\ell T_\ell,s_*\sigma\rangle,$$
with summation over $\ell\le n/2$. 
Let us summarize:
$$
\langle\sum_\epsilon(\epsilon s)^*T_k,\sigma\rangle=2^{N-n}{n+1\choose k}(-1)^k\sum_\ell\langle (-1)^\ell s^*T_\ell,\sigma\rangle.$$
It follows that
$$\sum_\epsilon(\epsilon s)^*T_k=2^{N-n}{n+1\choose k}(-1)^k\sum_\ell(-1)^\ell s^*T_\ell.$$
Now recall that, by Theorem 4.4, each $(\epsilon s)^*T_k$ is a cocycle representing
the cohomology class $\eu_k(P)$. Therefore
$$2^N\eu_k(P)=2^{N-n}{n+1\choose k}(-1)^k\sum_\ell(-1)^\ell \eu_\ell(P).$$
Comparing this formula for $k=0$ and for any other value of  $k$ we get the
the following Lemma
(which may be regarded as a variant of Theorem \a.1).
\smallskip
\bf Lemma \a.\t.\sl\par
Let $P$ be a flat principal $\PpGp$-bundle over a finite simplicial
complex $Y$ that has $N$ vertices. Then
$$2^N\eu_k(P)=2^N(-1)^k{n+1\choose k}\eu_0(P).$$\rm
\smallskip
Evaluating both sides of the formula from the Lemma on $[Z]$
concludes the proof of Theorem \a.1.\qed
\smallskip
\bf Corollary \a.\t.\sl\par
Let $P$ be a $P_+GL_+(n,K)$-bundle over an even-dimensional manifold $M^n$.
Then any triangulation of $M$ has at least $2^n|\langle \eu_0(P),[M]\rangle|$ simplices of dimension $n$.\rm
\smallskip
Proof. Pick a generic section $s$, over the given triangulation, of the associated bundle with fibre $\P_+$.
Then, by Theorem \a.1, $|\langle s^*T_k,[M]\rangle|=|\langle\eu_k(P),[M]\rangle|={n+1\choose k}|\langle\eu_0(P),[M]\rangle|$.
Since $\langle s^*T_k,\sigma\rangle=\pm1$ and the supports of the cocycles $s^*T_k$ are pairwise disjoint,
the number of $n$-simplices of the triangulation is at least
$$\sum_{k=0}^{\lfloor n/2\rfloor}|\langle s^*T_k,[M]\rangle|=
|\langle\eu_0(P),[M]\rangle|\cdot\sum_{k=0}^{\lfloor n/2\rfloor}{n+1\choose k}
=2^n|\langle\eu_0(P),[M]\rangle|.$$
\qed
\medskip
\medskip
\bf 11. Cross product of Euler classes.\rm
\medskip
\def\a{11}
\m=1
\n=1
\medskip
It will be convenient to put $[n]=\{0,1,\ldots,n\}$. We will use groups $GL_+(n,K)$ for varying $n$, therefore
we denote $U_+$ by $U_{n,+}$ in this and the next section.
\smallskip
\bf Theorem \a.\t.\sl\par
Let $E$ and $E'$ be $GL_+(n,K)$- and $GL_+(k,K)$-bundles over simplicial complexes $X$, $X'$ respectively.
Let $E\times E'$ be the product bundle over $X\times X'$. For any simplicial cycles
$z\in Z_n(X,\Z)$, $z'\in Z_k(X',\Z)$, we have
$$\langle \eu_0(E),z\rangle\cdot\langle \eu_0(E'),z'\rangle=\langle \eu_0(E\times E'),z\times z'\rangle.\leqno(\a.\f)$$
\par
\rm Proof. We first explain the general strategy of the proof.
We may and do assume that $X=\supp{z}$, $X'=\supp{z'}$.
We triangulate $X\times X'$ subdividing each product of simplices
$\sigma\times \sigma'$ in a standard way (to be recalled later).
It is also convenient to treat $E$, $E'$ and $E\times E'$ not as principal
bundles, but as vector bundles; e.g.~the fibre $E_x$ will be an $n$-dimensional
vector space over $K$. We pick generic sections: $s$ of $E$ and $s'$ of $E'$,
and combine them to a section $S$ of $E\times E'$. To ensure genericity
of $S$ we impose stronger than usual, weird-looking genericity conditions
on $s$ and on $s'$. The section $s$ induces a simplicial section
$s_+$ of the associated $X_{n,+}$-bundle, where $X_{n,+}$ is the
complex of generic tuples of points in $\P^{n-1}_+(K)$. Then, by Theorem 4.4,
we get cocycles $s^*_+T_+$ and $s^*_+T_0$ representing $\eu_+(E)$ and
$\eu_0(E)$. For each $n$-simplex $\sigma$ in $X$ we have
$\langle s^*_+T_+,\sigma\rangle=k^\pm$ for some $k$: if $k=0$, then
$\langle s^*_+T_0,\sigma\rangle=\pm1$; otherwise
$\langle s^*_+T_0,\sigma\rangle=0$. Similarly, we have cocycles
$s'^*_+T_0$ and $S^*_+T_0$ representing $\eu_0(E')$ and
$\eu_0(E\times E')$. Suppose that $z=\sum_\sigma n_\sigma\sigma$ and
$z=\sum_{\sigma'} n_{\sigma'}\sigma'$; then
$z\times z'=
\sum_{\sigma,\sigma'}n_\sigma n_{\sigma'}\cdot\sigma\times\sigma'$,
where $\sigma\times \sigma'$ denotes the chain representing the
standard subdivision of the product of simplices. We have
$$\langle \eu_0(E),z\rangle=
\sum_\sigma n_\sigma\langle s^*_+T_0,\sigma\rangle,\qquad
\langle \eu_0(E'),z'\rangle=\sum_{\sigma'}n_{\sigma'}\langle s'^*_+T_0,\sigma'
\rangle,$$
$$\langle \eu_0(E\times E'),z\times z'\rangle=
\sum_{\sigma,\sigma'} n_\sigma n_{\sigma'}\langle s^*_+T_0,
\sigma\times\sigma'\rangle.$$
Thus, to establish the theorem
it is enough to show that
$$\langle s^*_+T_0,\sigma\rangle\cdot\langle s'^*_+T_0,\sigma'\rangle
=\langle S^*_+T_0,\sigma\times\sigma'\rangle.\leqno(\a.\f)$$
We do this step-by-step. In Corollary \a.4 we show that if the left hand side
of (\a.2) is zero, then so is the right hand side. In Corollary \a.5 we show
that if the left hand side is non-zero, then in
the chain $\sigma\times\sigma'$ there is a unique summand
(unique $(n+k)$-simplex) on which $S^*_+T_0$ evaluates to $\pm1$.
Finally, in Lemma \a.6 we check that the sign of that evaluation is consistent
with (\a.2).

We proceed to the details.
First we pick a generic section $s$ of $E$ over $X^{(0)}$. The genericity condition is as follows.
For each $\ell$-simplex $\sigma=(x_0,\ldots,x_\ell)$ of $X$ ($\ell<n$),
the vectors $(s(x_0),\ldots,s(x_\ell))$ are linearly independent. For each
$n$-simplex $\sigma=(x_0,\ldots,x_n)$, if $\sum_{i=0}^n\alpha_is(x_i)=0$ is a non-trivial
linear relation (projectively unique,
because of the previous condition), then $\sum_{i=0}^n\alpha_i\ne0$.
To make sense of these conditions we choose a (flat)
trivialization of $E$ over $\sigma$.

This kind of section can be chosen vertex--by--vertex. Let $X^{(0)}=(x_1,x_2,\ldots,x_N)$.
First, we choose any non-zero $s(x_1)\in E_{x_1}$. When $s(x_1),\ldots,s(x_{i-1})$ have been chosen,
we choose (flat) trivializations over all simplices with vertex $x_i$. If $\sigma=(x_i,y_1,\ldots,y_\ell)$ is an
$\ell$-simplex of $X$ ($\ell<n$) such that $s(y_1),\ldots,s(y_\ell)$ have already been chosen, we use
the trivialization of $E$ over $\sigma$ to transport all $s(y_j)$ to $E_{x_i}$. There, these vectors
span a linear subspace $E_i^\sigma$ of dimension $\ell<n$. We have to ensure $s(x_i)\not\in E_i^\sigma$
in order to fulfill the first genericity condition for $\sigma$.

Next, for each simplex $\sigma=(x_i,y_1,\ldots,y_n)$ with $s(y_1),\ldots,s(y_n)$ already chosen
we pick a (flat) trivialization of $E$ over $\sigma$ and use it to transport the $s(y_j)$ to vectors
$s_j^\sigma\in E_{x_i}$. Then we form an affine subspace:
$$E_i^\sigma=\{\alpha_1s_1^\sigma+\ldots+\alpha_ns_n^\sigma\mid \alpha_1+\ldots+\alpha_n=1\}.$$
We have to choose $s(x_i)$ outside of this subspace in order to fulfill the second
genericity condition for $\sigma$.

A linear space over an ordered (hence infinite) field is not a union
of finitely many proper affine subspaces. Therefore, $s(x_i)$ can be suitably chosen.
By induction, there exists a generic section $s$ of $E$ over $X^{(0)}$.

With the section $s$ we associate a collection of scalars $\cal{A}$.
For each $n$-simplex $\sigma=(x_0,\ldots,x_n)$ of $X$ let $\sum_{i=0}^n\alpha_is(x_i)=0$
be the linear relation with $\sum_{i=0}^n\alpha_i=1$ (in some trivialization of $E$ over $\sigma$).
For every proper non-empty subset of $[n]$ we sum the corresponding $\alpha_i$'s.
The set $\cal{A}$ is the collection of all such sums (over all $n$-simplices).

Now, analogously, we choose a generic section $s'$ of $E'$ over $X'^{(0)}$.
It has its own collection of scalars $\cal{A'}$. We want $\cal{A}$ and $\cal{A'}$
to be disjoint. To this end, we perform the above section--choosing procedure for
$E'$ with supplementary restrictions. Suppose that we are at step $i$, choosing
$s'(x'_i)$. There is a collection of proper affine subspaces in $E'_{x_i'}$ that we need to avoid;
we now describe an additional finite collection, that will enforce our extra
``joint genericity'' condition. Let $\sigma'=(x'_i,y'_1,\ldots,y'_k)$ be a $k$-simplex of $X'$,
such that $s'(y'_j)$ are already chosen; let $s_j^{\sigma'}$ be $s'(y'_j)$ transported
to $E'_{x'_i}$ via a chosen trivialization of $E'$ over $\sigma'$.
For any generic (in the previous sense) $s'_i=s'(x'_i)$ there is a unique
relation $\beta_0s'_i+\sum_{j=1}^k\beta_js_j^{\sigma'}=0$ satisfying $\sum_{j=0}^k\beta_j=1$.
Pick an $\alpha\in\cal{A}$ and a proper non-empty $J\subset[k]$; we want to ensure
that $\sum_{j\in J}\beta_j\ne\alpha$. Let us express this as a restriction for the possible position
of $s'_i$. Suppose that $s'_i=\sum_{j=1}^k\gamma_js_j^{\sigma'}$, and that $\sum_{j\in J}\beta_j=\alpha$.
Let us express $\beta_j$ in terms of the $\gamma_j$. By the original genericity requirement
we know that $\Gamma:=-1+\sum_{j=1}^k\gamma_j\ne0$; therefore
$$-{1\over \Gamma}s'_j+\sum_{i=1}^k{\gamma_i\over\Gamma}s_j^{\sigma'}=0,$$
so that $\beta_0=-1/\Gamma$, $\beta_j=\gamma_j/\Gamma$. Thus, the condition $\sum_{j\in J}\beta_j=\alpha$
can be rewritten in terms of the $\gamma_j$ (putting $\gamma_0=-1$):
$\sum_{j\in J}\gamma_j/\Gamma=\alpha$, or
$\sum_{j\in J}\gamma_j=\alpha\Gamma$, or finally:
$$\sum_{j=0}^k(\alpha-\delta_J(j))\gamma_j=0.$$
Since $J$ is proper and non-empty, regardless of the value of $\alpha$ the set of vectors
$\sum_{j=1}^k\gamma_js_j^{\sigma'}$ for $\gamma_j$ satisfying this condition
forms a proper affine subspace of $E'_{x'_i}$. (The two special suspect cases: $J=\{0\}$, $\alpha=0$
and $J=\{1,\ldots,k\}$, $\alpha=1$ are easily seen to be impossible.)
Thus, the extra genericity conditions produce a new finite collection of proper affine subspaces
to avoid, so that it is possible to fulfill them.

Assume then that we have chosen jointly--generic (in the above sense) sections---$s$ of $E$ and $s'$ of $E'$.
We now form a generic section $S$ of $E\times E'$ over $(X\times X')^{(0)}$ by
$S(x,x')=(s(x),s'(x'))$. To claim genericity, we need to describe the (standard) triangulation
of $X\times X'$. We choose some total orders on $X^{(0)}$ and on $X'^{(0)}$, and order each simplex
of $X^{(n)}$ and of $X'^{(k)}$ accordingly. Let $\sigma=(x_0,\ldots,x_n)\in X^{(n)}$, and let
$\sigma'=(x'_0,\ldots,x'_k)\in X'^{(k)}$. Let $x_{(i,j)}=(x_i,x'_j)$. We form the $n\times k$ integer
grid---with vertex set $[n]\times[k]$ and edges connecting pairs that differ on exactly one coordinate
and exactly by 1. Shortest paths from $(0,0)$ to $(n,k)$ will be called \it admissible. \rm
(``Shortest'' is equivalent to ``going right or up at each step''.)
For each admissible path $\gamma\colon[n+k]\to[n]\times[k]$ we span an $(n+k)$-simplex $\sigma_\gamma$
in $\sigma\times\sigma'$ on the vertices $(S(x_{\gamma(j)})\mid j\in[n+k])$.
It is well-known that the set of all such $\sigma_\gamma$ triangulates $\sigma\times\sigma'$
(see [GM, I.1.5]).

We will call an $(n+1)$-tuple of vectors in an $n$-dimensional vector space
\it linearly generic, \rm if every $n$ of them are linearly independent.
\smallskip
\bf Lemma \a.\t.\sl\par
Vectors $(v_0,\ldots,v_n)$ are linearly generic if and only if the following
condition holds:
there is a projectively unique linear relation $\sum_{i=0}^n\alpha_iv_i=0$, and the 
coefficients in this relation are all non-zero.
\par
\rm Proof.
$(\Leftarrow)$ If some $n$ of the $v_i$'s were linearly dependent, a non-trivial linear
relation between them could be extended---by adding $0$ times the remaining vector---to
a non-trivial relation with coefficient $0$. Contradiction.

$(\Rightarrow)$ For dimensional reasons, there is a non-trivial linear relation
between the $v_i$'s; if some of its coefficients were $0$, it would give linear
dependence of a proper subset of the $v_i$'s. If the relation was not projectively
unique, one could form a linear combination of two non-proportional relations
and obtain a non-trivial relation with coefficient $0$.\qed{(Lemma)}
\smallskip
Observe that for a linearly generic tuple $(v_0,\ldots,v_n)$ the class $[(v_0,\ldots,v_n)]$ in $U_{n,+}$
is $0^\pm$ if and only if all the coefficients in the linear relation $\sum_{i=0}^n\alpha_iv_i=0$
are of the same sign.

Now we will tackle the question of genericity of the section $S$ (of $E\times E'$ over $(X\times X')^{(0)}$).
Let $\sigma=(x_0,\ldots,x_n)\in X^{(n)}$,  $\sigma'=(x'_0,\ldots,x'_k)\in X'^{(k)}$. Using trivializations
of $E$ over $\sigma$ and of $E'$ over $\sigma'$ we identify each $E_{x_i}$ with the same vector space
$V(\cong K^k)$, and each $E'_{x'_j}$ with $W(\cong K^k)$. Thus, we put $v_i=s(x_i)\in V$, $w_j=s'(x'_j)\in W$,
$V_{(i,j)}=(v_i,w_j)=S(x_{(i,j)})\in V\oplus W$. We would like to show that for each admissible path $\gamma$
the vectors $(V_{\gamma(j)}\mid j\in[n+k])$ are linearly generic. There are unique scalars $\alpha_i$ and $\beta_j$
such that
$$\sum_{i=0}^n\alpha_iv_i=0,\quad \sum_{i=0}^n\alpha_i=1;\qquad 
\sum_{j=0}^k\beta_jw_j=0,\quad \sum_{j=0}^k\beta_j=1.$$
Let $A_u=\sum_{i=0}^u\alpha_i$, $B_s=\sum_{j=0}^s\beta_j$. For a given path $\gamma$, let us arrange these
two sequences into one:
$$C_j=\cases{
A_i, & if $\gamma(j)=(i,*)$ and $\gamma(j+1)=(i+1,*)$;\cr
B_i, & if $\gamma(j)=(*,i)$ and $\gamma(j+1)=(*,i+1)$.
}$$
We put $C_{n+k}=A_n=B_k=1$ and $C_{-1}=0$.
\smallskip
\bf Lemma \a.\t.\sl\par
There is a projectively unique linear relation between the vectors
$(V_{\gamma(j)}\mid j\in[n+k])$. If we require that the sum of coefficients
be $1$, this relation
reads:
$$\sum_{j=0}^{n+k}(C_j-C_{j-1})V_{\gamma(j)}=0.$$
If all $A_u$ and $B_s$ are distinct, all coefficients of this relation are non-zero.
\par
\rm Proof. First, let us prove the formula. Projecting onto the first factor we get:
$$\sum_{i=0}^n(\sum_{j\in\Gamma_i}(C_j-C_{j-1}))v_i=0,$$
where $\Gamma_i=\{j\in[n+k]\mid \gamma(j)=(i,*)\}$ (similarly, we put
$\Gamma^i=\{j\in[n+k]\mid \gamma(j)=(*,i)\}$). We have $\Gamma_i=\{u,u+1,\ldots,u+\ell\}$ 
for some integers $u$, $\ell$. Therefore $\sum_{j\in\Gamma_i}(C_j-C_{j-1})=C_{u+\ell}-C_{u-1}=A_i-A_{i-1}=\alpha_i$.
Consequently, the displayed sum equals $\sum_{i=0}^n\alpha_iv_i=0$.
Similarly, the projection onto $W$ is $0$. Thus, the relation stated in the Lemma holds.

Suppose now that $\sum_{j=0}^{n+k}\Xi_jV_{\gamma(j)}=0$ is a linear relation with
$\sum_{j=0}^{n+k}\Xi_j=1$. Projecting onto the first factor, we get $\sum_{i=0}^n\xi_iv_i=0$,
where $\xi_i=\sum_{j\in\Gamma_i}\Xi_j$. Since $\sum_{i=0}^n\xi_i=\sum_{j=0}^{n+k}\Xi_j=1$ and
the vectors $v_i$ are linearly generic, we know that $\xi_i=\alpha_i$. Thus,
$$\sum_{j\in\Gamma_i}\Xi_j=\alpha_i.$$ Similarly, $$\sum_{j\in\Gamma^i}\Xi_j=\beta_i.$$
Since each $j$ is the largest element of exactly one set $\Gamma_i$ or $\Gamma^i$,
these equations recursively and uniquely determine all the $\Xi_j$.

Consequently, a linear relation between the $V_{\gamma(j)}$ with non-zero sum of coefficients
is projectively unique. This implies that there is no non-trivial relation with
sum of coefficients $0$---if it existed, it could be added to the one with sum of
coefficients $1$, contradicting the uniqueness of the latter.

The last claim of the lemma follows directly from the formula.\qed{(Lemma)}
\smallskip
\bf Corollary \a.\t.\sl\par
Suppose that the class $[(v_0,\ldots,v_n)]$ in $U_{n,+}$,
or the class $[(w_0,\ldots,w_k)]$ in $U_{k,+}$,
is not $0^\pm$. Then, for every admissible path $\gamma$, the class $[(V_{\gamma(j)}\mid j\in[n+k])]$
in $U_{n+k,+}$ is not $0^\pm$. \rm\par
Proof. The assumption can be interpreted as: $\alpha_i<0$ for some $i$, or $\beta_j<0$ for some $j$.
In each case one of the sequences $(A_u)$, $(B_s)$ is not increasing; therefore,
independently of $\gamma$, the sequence $(C_j)$ is not increasing. Consequently,
the relation between the $V_{\gamma(j)}$ (as in the Lemma) cannot have all positive coefficients,
while it does have some since their sum is $1$. Hence the claim.\qed
\smallskip
\bf Corollary \a.\t.\sl\par
Suppose that the class $[(v_0,\ldots,v_n)]=0^\pm$ in $U_{n,+}$, and the class $[(w_0,\ldots,w_k)]=0^\pm$ in $U_{k,+}$.
Then there is a unique admissible path $\gamma$ such that $[(V_{\gamma(j)}\mid j\in[n+k])]=0^\pm$ in $U_{n+k,+}$.
\rm\par
Proof. The sequences $(A_u)$, $(B_s)$ are increasing. There is a unique $\gamma$ such that $(C_j)$ is increasing
as well---then $[(V_{\gamma(j)}\mid j\in[n+k])]=0^\pm$. For other $\gamma$ we conclude as in the previous corollary.
\qed
\smallskip
Now we know that (\a.2) holds up to sign. To finish the proof it
remains to work out the relation between the signs that appear in the exponents
in Corollary \a.5, and to check that this relation is consistent with (\a.2).
The cycle $z\times z'$ contains a triangulated version of the product
$\sigma\times\sigma'$, for $\sigma\in\supp{z}$ and $\sigma'\in\supp(z')$, in the form
$\sum_\gamma\sgn(\gamma)\sigma_\gamma$. The summation is over all admissible $\gamma$. The sign
$\sgn(\gamma)$ equals $(-1)^{A(\gamma)}$, where $A(\gamma)$ is the area of the part of the grid that lies under
the image of $\gamma$. In particular, if $\gamma$ goes along the lower edge and the right-hand edge
of the grid, the sign is $+1$. If we change $\gamma$ by moving one $\gamma(j)$ to the opposite
vertex of a $1\times 1$ square---and get an admissible $\gamma'$---then $\sgn(\gamma')=-\sgn(\gamma)$.
\smallskip
\bf Lemma \a.\t.\sl\par
Suppose that $[(v_0,\ldots,v_n)]=0^s$, and that $[(w_0,\ldots,w_k)]=0^{s'}$.
Then 
$[(V_{\gamma(j)}\mid j\in[n+k])]=\sgn(\gamma)\cdot0^{ss'}$ for the $\gamma$ from the previous corollary.
\rm\par
Proof. We choose orientations of the bundles $E$, $E'$; we get induced orientations of
$V$, $W$ and $V\oplus W$. With respect to some positively oriented bases of $V$, $W$ we have:
$\sgn\det(v_0,\ldots,v_{n-1})=s$, $\sgn\det(w_0,\ldots,w_{k-1})=s'$.
We will show that for every admissible $\gamma$ the following sign formula holds:
$${\sgn\det}_B(V_{\gamma(j)}\mid j\in[n+k-1])=\sgn(\gamma),$$ where the determinant is calculated with respect
to the basis $$B=((v_0,0),\ldots,(v_{n-1},0),(0,w_0),\ldots,(0,w_{k-1})).$$ (This claim implies the lemma.)

First, for the $\gamma$ with $A(\gamma)=0$, the determinant is
$$\left|
\matrix{
1&0&\ldots&0&a_0&a_0&\ldots&a_0\cr
0&1&\ldots&0&a_1&a_1&\ldots&a_1\cr
&\vdots&&&\vdots&&&\vdots\cr
0&0&\ldots&1&a_{n-1}&a_{n-1}&\ldots&a_{n-1}\cr
1&1&\ldots&1&1&0&\ldots&0\cr
0&0&\ldots&0&0&1&\ldots&0\cr
&\vdots&&&\vdots&&&\vdots\cr
0&0&\ldots&0&0&0&\ldots&1}
\right|,
$$
where all $a_i$ are negative ($v_n=\sum_{i=0}^{n-1}a_iv_i$).
To calculate it, we use lower rows to cancel all the $a_i$ except the ones
in the $(n+1)^{\rm st}$ column. Then we use the left columns to cancel all
the remaining $a_i$---this increases the $(n+1,n+1)$-entry. The result is now
lower-triangular and positive on the diagonal.

Now let us consider the change of the determinant as $\gamma(j)$ moves across
a $1\times 1$ square. This changes one column. That column, and the neighbouring ones,
are as follows:
$$(\ldots,(v_i,w_j),(v_{i+1},w_j),(v_{i+1},w_{j+1}),\ldots)\leftrightarrow
(\ldots,(v_i,w_j),(v_i,w_{j+1}),(v_{i+1},w_{j+1}),\ldots).$$
The change, up to sign, can be performed by two column operations:
$$-(v_i,w_{j+1})=(v_{i+1},w_j)-(v_i,w_j)-(v_{i+1},w_{j+1}).$$
Since every admissible $\gamma$ can be obtained by such operations from the one
with $A(\gamma)=0$, the sign formula holds for all admissible paths.
\qed{(Lemma, Theorem \a.1)}
\medskip
\bf 12. Cup product of Euler classes.\rm
\medskip
\def\a{12}
\m=1
\n=1
Let $E$ be a (flat) $GL_+(n,K)$-bundle over a simplicial complex $X$.
We will often trivialize this bundle over simplices of $X$; to facilitate
the use of such trivializations we introduce the following convention.
Let $\sigma=(x_0,\ldots,x_\ell)$ be a simplex of $X$. We put $E_\sigma:=E_{x_0}$,
and we use any (flat) trivialization of $E$ over $\sigma$ to isomorphically
identify all the other $E_{x_i}$ with $E_\sigma$. Thus, if $s\colon X^{(0)}\to E$
is a section, we write $s(x_0),\ldots,s(x_\ell)\in E_\sigma$.
\smallskip
\bf Definition \a.\t. \rm A section $s\colon X^{(0)}\to E$ is called positive,
if for every simplex $\sigma=(x_0,\ldots,x_\ell)$ of $X$ there 
is a functional $\phi_\sigma\in E^*_\sigma$ such that $\phi_\sigma(s(x_i))>0$ for
$i=0,\ldots,\ell$.
\smallskip
If a $GL_+(n,K)$-bundle $E$ over $X$ admits a generic positive section $s$,
then $\langle \eu_0(E),z\rangle=0$ for every cycle $z\in Z_n(X,\Z)$.
Indeed, for every simplex $\sigma\in X^{(n)}$ we have
$s_*\sigma\ne 0^\pm$ in $U_{n,+}$, since the values of $s$ at the vertices of
$\sigma$ do not admit a linear relation with all positive coefficients---by
positivity of $s$. It turns out that (over a cycle) every positive section can be
perturbed to a generic positive section.
\smallskip
\bf Lemma \a.\t. \sl\par
If a $GL_+(n,K)$-bundle $E$ over a finite simplicial complex $X$ admits a positive section,
then it admits a generic positive section.
\rm \par
Proof. Let $s$ be a positive section, as witnessed by functionals $\phi_\sigma\in E^*_\sigma$
($\sigma\in X^{(n)}$). We will construct, vertex--by--vertex, a new generic section $s'$,
positive with respect to the same collection of functionals. We order the vertices of $X$,
and we start with $s'(x_0)=s(x_0)$. Suppose that $s'(x_\ell)$ have already been chosen for $\ell<i$.
Put $V=E_{x_i}$. When choosing $s'(x_i)$ in $V$, in order to ensure genericity,
we need to avoid a finite collection of affine hyperplanes, say defined by equations
$(\psi_j(v)=\alpha_j)_{j\in J}$
(where $\psi_j\in V^*$, $\alpha_j\in K$).
Also, for each $n$-simplex $\sigma$ with vertex $x_i$, we need to ensure that
$\phi_\sigma(s'(x_i))>0$ (we identify $E_\sigma$ with $V$). Let $w\in V$ be such that
$\psi_j(w)\ne 0$ for all $j\in J$; such $w$ exists, since $V$ is not the union of
finitely many hyperplanes $(\ker{\psi_j})_{j\in J}$. We will find suitable
$s'(x_i)$ in the form $v(\alpha):=s(x_i)+\alpha w$, for some scalar $\alpha$.
First, observe that the equation $\psi_j(v(\beta))=\alpha_j$ has a unique solution:
$\beta_j=(\alpha_j-\psi_j(s(x_i)))/\psi_j(w)$. Let $B:=\min\{\beta_j\mid\beta_j>0\}$.
The condition $\phi_\sigma(v(\beta))>0$, i.e.~$\phi_\sigma(s(x_i))+\beta\phi_\sigma(w)>0$,
is equivalent to
$\beta>-\phi_\sigma(s(x_i))/\phi_\sigma(w)$ (in case $\phi_\sigma(w)>0$) or to 
$\beta<-\phi_\sigma(s(x_i))/\phi_\sigma(w)$ (in case $\phi_\sigma(w)<0$). 
We know that $\beta=0$ satisfies all these inequalities. Therefore, the scalar
$$M:=\min\{-\phi_\sigma(s(x_i))/\phi_\sigma(w)\mid \phi_\sigma(w)<0\}$$ is positive.
We put $\alpha:={1\over4}\min(B,M)$ and $s'(x_i)=v(\alpha)$.\qed
\smallskip

\bf Corollary \a.\t. \sl\par
Let $E$ and $E'$ be $GL_+(n,K)$- and $GL_+(k,K)$-bundles over simplicial complexes $X$, $X'$ respectively.
For any simplicial cycles
$z\in Z_{n-\ell}(X,\Z)$, $z'\in Z_{k+\ell}(X',\Z)$, where $\ell>0$, we have
$$\langle \eu_0(E\times E'),z\times z'\rangle=0.$$
\rm\par
Proof. We may and do assume that $X=\supp{z}$, $X'=\supp{z'}$.
Let $s$ be a generic section of $E$. For dimensional reasons, the values of
$s$ at the vertices of any simplex $\sigma$ of $X$ are linearly independent;
therefore, a functional $\phi_\sigma$ can be chosen that evaluates to $1$ on each of them.
Thus, $s$ is positive.
Now define $S\colon (X\times X')^{(0)}\to E\times E'$ by $S(x,x')=(s(x),0)$.
Then, for any simplices $\sigma\in X^{(n-\ell)}$ and $\sigma'\in X'^{(k+\ell)}$,
 and any $(n+k)$-dimensional simplex $\sigma_\gamma$ in the standard triangulation
 of $\sigma\times \sigma'$, we may put $\phi_{\sigma_\gamma}=\phi_\sigma\circ\pi_E$. Then,
for every vertex $(x,x')$ of $\sigma_\gamma$ we have
$$\phi_{\sigma_\gamma}(S(x,x'))=\phi_\sigma(\pi_E(s(x),0))=\phi_\sigma(s(x))>0.$$
Therefore, $S$ is a positive section of $E\times E'$ over $\supp{z\times z'}$.
By the lemma above, there exists a generic positive section, and that implies
the asserted vanishing.\qed
\smallskip
\bf Corollary \a.\t.\sl\par
Let $E$ and $E'$ be $GL_+(n,K)$- and $GL_+(k,K)$-bundles over simplicial complexes $X$, $X'$ respectively.
For any simplicial cycle
$Z\in Z_{n+k}(X\times X',\Z)$ we have
$$\langle \eu_0(E)\times \eu_0(E'),Z\rangle=\langle \eu_0(E\times E'),Z\rangle.$$
\rm\par
Proof. Indeed, by K\"unneth's formula, an integer multiple of $Z$ is 
homologous to a combination  of cycles of the form $z\times z'$; for the latter,
the formula holds either by the previous corollary, or by the theorem from the
previous section.\qed
\smallskip
\bf Theorem \a.\t.\sl\par\nobreak
Let $E$ and $E'$ be $GL_+(n,K)$- and $GL_+(k,K)$-bundles over a simplicial complex $X$.
For any simplicial cycle $z\in Z_{n+k}(X,\Z)$ we have
$$\langle \eu_0(E)\cup \eu_0(E'),z\rangle=\langle \eu_0(E\oplus E'),z\rangle.$$
\rm\par
Proof. Let $\Delta\colon X\to X\times X$ be the diagonal map.
$$\eqalign{
\langle \eu_0(E)\cup \eu_0(E'),z\rangle&=
\langle \Delta^*(\eu_0(E)\times \eu_0(E')),[z]\rangle\cr
&=\langle \eu_0(E)\times \eu_0(E'),\Delta_*[z]\rangle=
\langle \eu_0(E\times E'),\Delta_*[z]\rangle\cr
&=\langle \Delta^*\eu_0(E\times E'),[z]\rangle=
\langle \eu_0(\Delta^*(E\times E')),[z]\rangle\cr
&=\langle \eu_0(E\oplus E'),z\rangle.}$$
\qed
\medskip
\bf 13. Comparison of Euler and Witt classes.\rm
\medskip
\def\a{13}
\m=1
\n=1
\medskip
We use the functoriality theorem (Theorem 1.5) to compare various tautological classes
that we have constructed. We begin with $\eu$ and $\eu_+$.
\smallskip
\bf Euler classes. \rm We assume $n$ even. There is a natural map $\P_+\to\P$; it induces a
simplicial (non-degenerate) map $f\colon X_+\to X$. The groups $\PpGp$ and $\PGp$ acting on $X_+$ and $X$ (respectively)
are also related by the natural projection homomorphism
$\phi\colon\PpGp\!\to\!\PGp$.
The map $f$ is $\phi$-equivariant,
and induces a coefficient group map $f\colon U_+\to U$. Theorem 1.5 applies and gives the following diagram:
$$
H^n(\PGp,U)
\mathop{\longrightarrow}\limits^{\phi^*}
H^n(\PpGp,U)
\mathop{\longleftarrow}\limits^{f_*}
H^n(\PpGp,U_+).
$$ 
Recall that $U\simeq\Z$ and $U_+\simeq \Z^{(n/2)+1}$. The map $f\colon U_+\to U$
can be described explicitly using Remark 8.8. The generator $a^+$ of $U_+$ is represented by
the simplex $([e_1],\ldots,[e_n],[v_a])$, where $v_a=e_1+\ldots+e_a-(e_{a+1}+\ldots+e_n)$.
The image of this simplex in $X$ determines in $U$ the symbol $[\sgn(\det(e_1,\ldots,e_n)\cdot (-1)^{n-a}]=[(-1)^a]$.
Therefore $f(a^+)=[(-1)^a]=(-1)^a[+]$. It follows that the induced map on cohomology, $f_*\colon H^n(\PpGp,U_+)\to H^n(\PGp,U)$,
maps $\eu_+=\oplus_a\eu_a$ to $\sum_a(-1)^a\eu_a$. Theorem 1.5 implies the following result.
\smallskip
\bf Theorem \a.\t. \sl\par
Let $\phi\colon \PpGp\to\PGp$ be the natural projection homomorphism. Then
$$\phi^*\eu=\sum_{a=0}^{n/2}(-1)^a\eu_a.$$
\rm\par
A (flat) $\PpGp$-bundle $P$ over $Y$ determines a $\PGp$-bundle $P'$ over $Y$.
As is usual in such cases, we put $\eu(P):=\eu(P')\in H^n(Y,\Z)$.
\smallskip
\bf Corollary \a.\t.\sl\par
Let $P$ be a (flat) $P_+GL_+(n,K)$-bundle over an oriented closed $n$-manifold $M$.
Then
$$\langle\eu(P),[M]\rangle=2^n\langle\eu_0(P),[M]\rangle.$$
\rm\par
Proof. Using Theorem \a.1 and Theorem 10.1 we calculate
$$\eqalign{
\langle\eu(P),[M]\rangle&=
\langle\sum_{k=0}^{n/2}(-1)^k\eu_k(P),[M]\rangle=
\sum_{k=0}^{n/2}(-1)^k\langle\eu_k(P),[M]\rangle\cr
&=
\sum_{k=0}^{n/2}(-1)^k(-1)^k{n+1\choose k}\langle\eu_0(P),[M]\rangle=
2^n\langle\eu_0(P),[M]\rangle.}$$
\qed

\bf Remark \a.\t.\rm\par
For $n=2$ Theorem 9.1 gives $3\,\eu_0+\eu_1=0$. Theorem \a.1 now implies
$\phi^*\eu=4\,\eu_0$, i.e.~in this case Corollary \a.2 can be strengthened to equality
in $H^2(P_+GL_+(2,K),\Z)$---there is no need to evaluate on cycles.

\smallskip
\bf Witt class. \rm 
In Section 7 we discussed the action of $PSL(2,K)$ on $\P^1$, on the associated complex $X$,
and the resulting Witt class $w\in H^2(PSL(2,K),W(K))$. In Section 8 we considered
the action of $PGL_+(2,K)$ on the same spaces, and the resulting cohomology
class $\eu\in H^2(PGL_+(2,K),\Z)$. Theorem 1.5 may be applied to the identity map $\iota\colon X\to X$
and the injection homomorphism $\phi\colon PSL(2,K)\to PGL_+(2,K)$. Before stating the result we
compute the coefficient map $\iota\colon W(K)\to\Z$. The symbol $[\lambda]$ is represented by
the triple
$t_\lambda=\left(
{1\broose0},
{0\broose1},
{1\broose\lambda}
\right)$. To find the symbol
of $t_\lambda$ in $U_2(X,PGL_+(2,K))$ we write
${1\choose\lambda}=1\cdot{1\choose0}+\lambda\cdot{0\choose1}$;
then, using Remark 8.8, we get
$$\left[\sgn(\left|\matrix{1&0\cr0&1}\right|\cdot1\cdot\lambda)\right]=
[\sgn(\lambda)].$$
Therefore, the map $\iota\colon W(K)\to\Z$ is just the signature map $\sigma$,
given by $\sigma([\lambda])=\sgn(\lambda)$. The diagram and the theorem are as follows.
$$
H^2(PGL_+(2,K),\Z)
\mathop{\longrightarrow}\limits^{\phi^*}
H^2(PSL(2,K),\Z)
\mathop{\longleftarrow}\limits^{\sigma_*}
H^2(PSL(2,K),W(K))
$$
\bf Theorem \a.\t.\sl\par
Let $\phi\colon PSL(2,K)\to PGL_+(2,K)$ be the standard inclusion. Then
$$\phi^*\eu=\sigma_*w.$$
Furthermore, the pull-back of this class to $SL(2,K)$ is equal to $4\,\eu_0$.
\rm
\smallskip
The last claim of the theorem follows from Remark \a.3.
\smallskip
\bf Non-vanishing. \rm
Consider a flat vector $SL(2,\R)$-bundle $E$ over a closed oriented surface $\Sigma$.
The (classical, topological) Euler class $\eu_t(E)$ of $E$ (more precisely, the Euler number $\langle\eu_t(E),[\Sigma]\rangle$)
can be computed as the signed number of zeroes
of a generic section of $E$; generic means: transversal to the zero section.
Consider now a triangulation $Y$ of $\Sigma$. Let $s\colon Y^{(0)}\to E$ be a generic section
over the set of vertices of $Y$. Here genericity means that for every 2-simplex $\sigma$ of $Y$
the values of $s$ at the vertices of $\sigma$ are pairwise linearly independent
(as usual, we compare them using a flat trivialization of $E$ over $\sigma$).
The section $s$ can be affinely extended to each simplex of $Y$. Together, these extensions
define a generic section of $E$ over $\Sigma$ in the previous, classical sense.
Moreover, the zeroes of this extended section occur exactly in simplices $\sigma$ on which
$s^*T_0$ (the cocycle representing $\eu_0(E)$, see Remark 8.8 and Definition 8.9) is non-zero, 
and the sign of the zero in $\sigma$ is equal to $s^*T_0(\sigma)$. These arguments prove
the following statement.
\smallskip
\bf Fact \a.\t.\sl\par
Let $E$ be a flat $SL(2,\R)$-bundle over a closed surface $\Sigma$. Then
$$\langle\eu_0(E),[\Sigma]\rangle=\langle\eu_t(E),[\Sigma]\rangle.$$\rm
\smallskip
We will now prove that all
the Euler classes constructed in this paper are non-zero (for $n$ even).
\smallskip
\bf Theorem \a.\t.\sl\par
Let $K$ be an ordered field and let $n$ be even. Then the
Euler classes $\eu$, $\eu_+$ and all $\eu_k$ are non-zero.\rm
\smallskip
Proof. Recall that an ordered field contains $\Q$ as a subfield, and the order restricted to $\Q$ is standard.
Due to field restriction stability of our classes (see Remark 8.10) it is enough to show the theorem for $K=\Q$.

Assume first that $n=2$.
Recall that over a closed oriented surface $\Sigma$ of genus $g\ge2$
there are flat vector $SL(2,\R)$-bundles $E$ with non-trivial Euler class $\eu_t$
(see [MS, Appendix C]).
Moreover, Takeuchi proved that $SL(2,\Q)$ can be used as the structure
group of such bundles (see [Takeuchi]); let us call such examples (flat $SL(2,\Q)$-bundles with non-trivial $\eu_t$)  
Takeuchi bundles.
Fact \a.5 implies that the Euler class $\eu_0$ is non-zero for Takeuchi bundles.
Theorem 10.1 and Corollary \a.2 imply that also $\eu_1$ and $\eu$ are non-trivial for them.

For larger even $n=2k$ we consider the Cartesian product $Y$ of $k$ copies of $\Sigma$,
and over $Y$ the product bundle $E^{\times k}$ of $k$-copies of a Takeuchi bundle $E$.
Then Theorem 11.1 shows that $\langle \eu_0(E^{\times k}),[Y]\rangle=\langle \eu_0(E),[\Sigma]\rangle^k\ne0$.
Again, it follows from Theorem 10.1 and from Corollary \a.2 that all $\eu_k$ as well as $\eu$ are non-trivial
on $E^{\times k}$.\qed

\bigskip
\vfill\eject
\bf References\rm\nobreak
\medskip
\item{[BFG]}
 N.~Bergeron, E.~Falbel, A.~Guilloux,
 Tetrahedra of flags, volume and homology of SL(3). 
 Geom.~Topol. {\bf 18} (2014), no.~4, 1911--1971.

\smallskip
\item{[Brown]}
 K.~Brown,
 Cohomology of groups.
 Graduate Texts in Mathematics, 87.
 Springer-Verlag, New York-Berlin, 1982.

\smallskip
\item{[Dug]}
 J.~Dugundji,
 Cohomology of equivariant maps.
 Trans.~Amer.~Math.~Soc. {\bf 89} 1958, 408--420.

\smallskip
\item{[EKM]}
 R.~Elman, N.~Karpenko and A.~Merkurjev,
 The algebraic and geometric theory of quadratic forms.
 American Mathematical Society Colloquium Publications, 56.
 American Mathematical Society, Providence, RI, 2008. 

\smallskip
\item{[GM]} S.I.~Gelfand and Y.I.~Manin,
Methods of homological algebra. Second edition.
Sprin\-ger Monographs in Mathematics. Springer-Verlag, Berlin, 2003.

\item{[Ghys]} \'E.~Ghys, Groupes d'hom\'eomorphismes du cercle et
cohomologie born\'ee.
 In: The Lefschetz Centennial Conference, Part III (Mexico City, 1984),
 pp.~81--106. Amer.~Math.~Soc., Providence, RI, 1987.

\smallskip
\item{[Gro]}
 M.~Gromov,
 Volume and bounded cohomology.
 Inst.~Hautes \'{E}tudes Sci.~Publ.~Math., {\bf 56} (1982), 5--99.

\smallskip
\item{[Kr-T]}
 L.~Kramer, K.~Tent, A Maslov cocycle for unitary groups.
 Proc.~Lond.~Math.~Soc. (3) {\bf 100} (2010), no.~1, 91–115.

\smallskip
\item{[L\"oh]}
 C.~L\"oh,
 Group Cohomology \& Bounded Cohomology, 2010, preliminary version of lecture notes,
 
 {\tt http://www.mathematik.uni-regensburg.de/loeh/teaching/topologie3\_ws0910/prelim.pdf}

\smallskip
\item{[MS]}
 J.~Milnor, J.~Stasheff,
 Characteristic classes.
 Annals of Mathematics Studies, No.~76.
 Princeton University Press, Princeton, N.J.;
 University of Tokyo Press, Tokyo, 1974.

\smallskip
\item{[Morita]}
 S.~Morita,
 Geometry of Characteristic Classes.
 Translated from the 1999 Japanese original.
 Translations of Mathematical Monographs, 199.
 Iwanami Series in Modern Mathematics.
 American Mathematical Society, Providence, RI, 2001.

\smallskip
\item{[Ne]}
 J.~Nekov\'a$\check{\rm r}$,
 The Maslov index and Clifford algebras.
 (Russian) Funktsional.~Anal.~i Prilozhen. {\bf 24} (1990), no.~3, 36--44, 96;
 translation in Funct.~Anal.~Appl. {\bf 24} (1990), no.~3, 196--204 (1991). 

\smallskip
\item{[Rez]}
 A.~Reznikov,
 Euler class and free generation.
 arXiv:dg-ga/9709009.

\smallskip
\item{[Takeuchi]}
 K.~Takeuchi,
 Fuchsian groups contained in $SL_{2}(\Q)$,
 J.~Math.~Soc.~Japan, {\bf 23} (1971), no.~1, 82--94.

\bye